\documentclass[11pt,reqno]{amsart}
\usepackage[utf8]{inputenc}
\usepackage[numbers,sort]{natbib}
\usepackage[french,english]{babel}
\usepackage{amsfonts}       
\usepackage{nicefrac}       
\usepackage{microtype}      
\usepackage{lipsum}
\usepackage{graphicx}
\usepackage{tikz}
\usepackage{amsmath}
\usepackage{amsthm}
\usepackage{bbold}
\usepackage{algorithm}
\usepackage{algorithmic}
\usepackage{mathrsfs}
\usepackage{mathtools}
\usepackage{amssymb}
\usepackage{icomma}
\usepackage{comment}
\usepackage{caption}
\usepackage{subcaption}

\usepackage{pgfplots}
\pgfplotsset{compat=1.18}

\usetikzlibrary{patterns.meta}

\usepackage[colorlinks=true, pdfstartview=FitV, pagebackref=true]{hyperref}

\usepackage[a4paper,vmargin=3.5cm,inner=3cm,outer=3cm]{geometry}
\usepackage[dvipsnames,svgnames]{xcolor}

\usepackage{blkarray}
\newcommand{\lab}[1]{\scriptstyle #1}

\newtheorem{theorem}{Theorem}

\newtheorem{lemma}[theorem]{Lemma}
\newtheorem{prop}[theorem]{Proposition}

\newtheorem{definition}{Definition}

\newtheorem{remark}{Remark}

\newtheorem{assumption}{Assumption}

\usepackage{microtype}
\usepackage{graphicx}
\usepackage{wrapfig}
\usepackage{booktabs} 

\newcommand{\N}{\mathbb{N}}
\newcommand{\R}{\mathbb{R}}
\newcommand{\PP}{\mathbb{P}}

\renewcommand{\P}{\mathbb{P}}
\newcommand{\E}{\mathbb{E}}
\newcommand{\given}{\, |\, }
\newcommand{\cF}{\mathcal{F}}
\newcommand{\cM}{\mathcal{M}}
\newcommand{\cP}{\mathcal{P}}

 \newcommand{\ind}{\mathbb{1}}
 \newcommand{\e}{\operatorname{e}}
  \renewcommand{\d}{\mathrm{d}}
\newcommand{\model}{\textsc{Model}}
\renewcommand{\tilde}{\widetilde}

\usepackage{array,multirow,makecell}
\setcellgapes{1pt}
\makegapedcells
\newcolumntype{R}[1]{>{\raggedleft\arraybackslash }b{#1}}
\newcolumntype{L}[1]{>{\raggedright\arraybackslash }b{#1}}
\newcolumntype{C}[1]{>{\centering\arraybackslash }b{#1}}

\title{Sublinear multi-type preferential attachment}

\author{Anna BEN-HAMOU}
\address{Sorbonne Université, Université Paris Cité, CNRS, Laboratoire de Probabilités, Statistique et Modélisation, LPSM, F-75005 Paris, France}
\email{abenhamou@lpsm.paris}
\author{Nicolas BROUTIN}
\email{broutinn@lpsm.paris}
\author{Camille CAZAUX}
\email{ccazaux@lpsm.paris}

\begin{document}

\maketitle

\begin{abstract}
    We consider a random growing graph model with two \emph{types}, or \emph{communities}, that belongs to the class of preferential attachment models. Starting from an initial graph $G_0$ in which each node has a type, the graph $G_{n+1}$ is generated by addition of a new node, which attaches through an edge to a node $v$ of $G_n$. The probability of attaching to $v$ depends on the communities of both the new node and $v$, and is proportional to $\deg_n(v)^\gamma$, where $\gamma\in[0, 1]$. We establish a convergence result for the limiting proportion of vertices of each degree and in each community in a certain regime of the parameter $\gamma$. To this end, we first use stochastic approximation techniques to reduce the convergence problem to the study of a nonlinear evolution equation. We then use a probabilistic coupling approach to investigate the long-time behavior of the unique solution to this evolution equation. We emphasize that the techniques previously used to handle the single-community case cannot be extended directly to our setting. Indeed, the associated continuous-time embedding is not a multitype branching process in the usual sense, since the transition rates at any given time $t$ depend on the entire graph structure at that time.
\end{abstract}

\section{Introduction}

\subsection{Setting and main results}

We study a random growing graph model with two \emph{types}, or \emph{communities} ($A$ and $B$), that belongs to the class of preferential attachment models. Starting from an initial graph $G_0$ in which each node has a type, the graph $G_{n+1}$ is obtained by adding to $G_n$ a new node $v_{n+1}$. This node first receives a type drawn according to a Bernoulli distribution. Then, given its type, it attaches by an edge to a randomly chosen node $v$ --- belonging to the node set $V(G_n)$ of $G_n$ --- with a probability that depends on the type of~$v$ and sublinearly on its degree. More precisely, we consider the following random graph model.

\begin{definition}[\model\ $1$]\label{def:model-1}
Let $p, q_A, q_B$ and $\gamma$ be four parameters in $[0, 1]$. Start with a fixed finite connected graph $G_0$, whose $n_0 := \# V(G_0)$ nodes are of type $A$ or $B$. Then, for $n\geq 0$, the graph $G_{n+1}$ is obtained from $G_n$ by adding a new node $v_{n+1}$ endowed with an edge. The attachment rule is the following: 
\begin{enumerate}
\item[(i)] the type $\sigma(v_{n+1})$ of $v_{n+1}$ is randomly chosen: type $A$ with probability $p$, type $B$ with probability $1-p$;
\item[(ii)]
\begin{enumerate}
\item[(a)] given its own community, the new node chooses whether it will attach to a node in the same community with probability~$1-q_{\sigma(v_{n+1})}$, or in the other community with probability~$q_{\sigma(v_{n+1})}$;
\item[(b)] it attaches to a node $u$ in the community chosen at step~(a) with probability proportional to $\deg_n(u)^\gamma$, where $\deg_n(u)$ denotes the degree of $u$ in $G_n$. 
\end{enumerate}
\end{enumerate}
\end{definition}

One may also consider a ``softer'' variant, referred to as \model\ $2$, where step~(ii) is replaced by:
\begin{enumerate}
\item[(ii')] the new node attaches to a node $u$ of $G_n$ with probability proportional to 
\[
\begin{cases}
(1-q_{\sigma(v_{n+1})})\deg_n(u)^\gamma & \text{ if $\sigma(v_{n+1})=\sigma(u)$,}\\
q_{\sigma(v_{n+1})} \deg_n(u)^\gamma & \text{ if $\sigma(v_{n+1})\neq \sigma(u)$.}
\end{cases}
\]
\end{enumerate}

More precisely, to make the difference between \model\ $1$ and \model\ $2$ more explicit, let $\Lambda_{n+1} := \sigma(v_{n+1})$ denote the type of the new node arriving at time $n+1$. The random variables $(\Lambda_n)_{n\geq 1}$, are \emph{i.i.d.}\ with distribution $\mathrm{Bernoulli}(p)$. We then define the filtration $(\mathcal{F}_n)_{n \in \N}$ by setting $\mathcal{F}_n := \sigma((G_i)_{0 \leq i \leq n}, (\Lambda_n)_{1 \leq i \leq n})$, for each $n \in \N$. Let us also define, for $\Lambda\in\{A, B\}$, the random variable
\[
D_n(\Lambda) := \displaystyle\sum_{\substack{v \in V(G_n) \\ \sigma(v) = \Lambda}} \deg_n(v)^\gamma \, ,
\]
and for $\Lambda\in\{A, B\}$, let $\overline{\Lambda}$ denote the ``other'' type. In \model\ $1$, the conditional probability that the new node $v_{n+1}$ attaches to a given node $u \in V(G_n)$, given $\cF_n$ and $\Lambda_{n+1}$ is
\[
(1-q_{\Lambda_{n+1}})\frac{\deg_n(u)^\gamma}{D_n\left(\Lambda_{n+1}\right)}\ind_{\{\sigma(u)=\Lambda_{n+1}\}}
+q_{\Lambda_{n+1}}\frac{\deg_n(u)^\gamma}{D_n\left(\overline{\Lambda}_{n+1}\right)}\ind_{\{\sigma(u)\neq \Lambda_{n+1}\}}\, ,
\]
whereas in \model\ $2$, this probability is
\[
\frac{(1-q_{\Lambda_{n+1}})\deg_n(u)^\gamma\ind_{\{\sigma(u)=\Lambda_{n+1}\}}+q_{\Lambda_{n+1}}\deg_n(u)^\gamma\ind_{\{\sigma(u)\neq \Lambda_{n+1}\}}}{(1-q_{\Lambda_{n+1}})D_n\left(\Lambda_{n+1}\right)+q_{\Lambda_{n+1}}D_n\left(\overline{\Lambda}_{n+1}\right)}\, \cdot
\]

\begin{remark}
In both models, the parameters $q_A$ and $q_B$ govern the inclination to attach to the ``other'' community. Let us mention here that we will sometimes write $p_A$ and $p_B$ instead of $p$ and $1-p$. Finally, observe that one may recover the one-community case by setting $p=1$ and $q_A=0$.
\end{remark}

Our goal is to determine the limiting degree distribution in~$G_n$ when $n$ tends to~$+\infty$. We are thus led to consider the sequence of counting measures $(m_n)_{n\in\N}$ on $E := \N^* \times \{A, B\}$ given by\footnote{Throughout this paper, we denote by $\N^*=\{1, 2, 3, \dots\}$ the set of positive integers and by $\N=\{0\}\cup\N^*$ the set of nonnegative integers.}
\begin{align}\label{def-m-n}
\forall n\geq 0,\, \forall x \in E,\quad m_n(x) = \#\{v \in V(G_n) : (\deg_n(v), \sigma(v))= x\}\, .
\end{align}
In other words, if $x = (k, \Lambda)$, then $m_n(x)$ counts the number of nodes in~$G_n$ with degree $k$ and type $\Lambda$ and $\frac{m_n(x)}{n+n_0}$ denotes the proportion of such nodes in $G_n$. The counting measure $m_n$ will be identified to a sequence indexed by the lexicographic order on $E$: $(1, A), (1, B), (2, A), (2, B)...$

Our main contribution is to establish almost sure convergence in total variation distance of the renormalized counting measure $\frac{m_n}{n+n_0}$ (as $n$ goes to infinity): in a certain parameter regime $\gamma \in [0, \gamma_0)$ ($\gamma_0$ being a strictly positive constant in $(0, 1]$) in the two-community case (Theorem~\ref{th:deg-distrib-2}), and for any choice of $\gamma \in [0, 1]$ satisfying a rather remarkable inequality in the one-community case (Theorem~\ref{th:deg-distrib-1}). 

\begin{theorem}[Limiting degree distribution in the two-community case]\label{th:deg-distrib-2} 
Let $(G_n)_{n \geq 0}$ be the sequence of random graphs generated according to \model\ $i\in\{1, 2\}$ with parameters $p, q_A, q_B$ and $\gamma$ all in $[0, 1]$. There exists a probability measure $\mu^\star_{(i)}$ on $E$ and a constant $0 < \gamma_0 \leq 1$ such that for $\gamma\in [0, \gamma_0)$, we have
\[
\left\|\frac{m_n}{n_0+n}-\mu^\star_{(i)} \right\|_{\mathrm{TV}}\, \xrightarrow[n\to\infty]{\mathrm{a.s.}} \, 0\, .
\]
\end{theorem}

We denote by $\|\cdot\|_{\mathrm{TV}}$ the total variation distance, with the convention that $\|\cdot\|_{\mathrm{TV}}=\frac{1}{2}\|\cdot\|_1$, where $\|\cdot\|_1$ denotes the $\ell^1$-norm.

\begin{remark}
By the results of \citet{jordan2013geometric} and~\citet{bhamidi2022community}, this was already known for $\gamma=0$ and $\gamma=1$ (they showed pointwise convergence to a limiting probability distribution, which is equivalent to convergence in total variation by Scheffé's lemma, since the state space is countable). We conjecture that the result of Theorem~\ref{th:deg-distrib-2} holds for all $\gamma \in [0, 1]$. The precise value of our $\gamma_0$ is implicitly defined in~\eqref{condition-gamma-convergence}. Let us point out here that for $p=\tfrac{1}{2}$ and $q_A=q_B=q$, one can take $\gamma_0=\tfrac{1}{7}$. 
\end{remark}

\begin{theorem}[Limiting degree distribution in the one-community case]\label{th:deg-distrib-1}  
Let $(G_n)_{n \in \N}$ be the sequence of random graphs generated according to Definition~\ref{def:model-1} with parameters $p = 1$, $q_A = 0$ and $\gamma \in [0, 1]$, corresponding to the single-community case. Let $E = \N^*$ and, for $n \in \N$, define the counting measure $m_n$ by
\[
m_n(x) = \#\{v \in V(G_n) : \deg_n(v)= x\}\, ,
\]
for each $x \in E$. There exists a probability measure $\mu^\star$ on $E$ such that, for each $\gamma \in [0, 1]$ satisfying the following inequality
\begin{equation}\label{eq:inequality}
\sum_{k\geq 1} \prod_{j=1}^k\frac{j^\gamma}{j^\gamma +3^\gamma-1} > 1\,,
\end{equation}
we have
\[
\left\|\frac{m_n}{n_0+n}-\mu^\star \right\|_{\mathrm{TV}}\, \xrightarrow[n\to\infty]{\mathrm{a.s.}} \, 0\, .
\]
\end{theorem}

\subsection{Related work}

Complex networks arise in a wide variety of settings, including technological, biological, social, and information systems. In many real-world networks, empirical degree distributions exhibit a power-law tail, that is, the probability that a randomly chosen node has degree $k$ is approximately proportional to $k^{-\alpha}$, where $\alpha$ often lies between $2$ and $3$. Such networks typically contain many vertices with relatively few connections and a small fraction of highly connected vertices, known as hubs. This observation has motivated the study of growing network models, in which new vertices and edges are added over time. In 1999, Barabási and Albert introduced the \emph{preferential attachment} mechanism, whereby new vertices are more likely to connect to vertices that already have high degrees (see~\cite{barabasi1999emergence}). This self-reinforcing mechanism can give rise to power-law decay of the degree distribution, a precise mathematical formulation and proof of this result were given in~\cite{bollobas2001degree} and independently in~\cite{mori2002random}. The classical Barabási–Albert model is recovered as a special case of our model in the single-community setting when $\gamma=1$. This model is now well understood. By contrast, in the classical Erdős–Rényi model $G(n, p)$, node degrees follow a binomial distribution and, in the usual regimes, concentrate around a characteristic scale. Power-law degree distributions, by comparison, are heavy-tailed and are commonly described as \emph{scale-free}, reflecting the absence of a characteristic degree scale.

Several extensions of the model have been studied. The case of uniform attachment, corresponding to $\gamma=0$ in our setting, is well understood and corresponds to the uniform recursive tree. When $\gamma>1$, a condensation phenomenon occurs: in the limit, a single random node dominates the network by accumulating a disproportionate number of connections. This phenomenon was already known in the physics literature, notably through the work of Krapivsky, Redner, and Leyvraz (see~\cite{krapivsky2000connectivity}), before being established rigorously in subsequent mathematical studies (see, for example,~\cite{oliveira2005connectivity}, and more recently~\cite{sethuraman2019growth}). More generally, the probability governing the choice of the node or vertices to which a new node attaches may depend on the degrees in a variety of ways. Given a weight function $w:\N \to \R_+$, the attachment probability may be proportional to $w$ evaluated at the degree of a node. In~\cite{oliveira2005connectivity} and~\cite{sethuraman2019growth}, the weight function $w$ is superlinear, and the techniques developed therein rely on martingale arguments that are particularly well suited to this setting. The paper~\cite{athreya2007preferential} studies preferential attachment graphs with specific sublinear weight functions $w(k)=ck^\gamma$ (with $\gamma \in (\frac{1}{2}, 1)$), linear functions $w(k)=ak+b$, and superlinear functions, showing that they lead to distinct asymptotic behaviors: slow degree growth, power-law degree distributions, and extreme concentration of connections, respectively. These results are established using probabilistic methods, particularly continuous-time embeddings and pure birth processes.

Using probabilistic asymptotic techniques based on embeddings of random tree growth processes into continuous-time general branching processes, \citet{rudas2004random} and~\citet{rudas2007random} determine the limiting degree distribution and the local limit in the sublinear regime.  In~\cite{rudas2007random}, the main results rely on branching process limit theorems and random characteristics to establish the almost-sure convergence of the empirical distributions of node degrees and descendant subtrees. By contrast, to the best of our knowledge, the limiting degree distribution of the model in the sublinear regime with communities has not yet been investigated.

Meanwhile, when $\gamma = 0$ or $\gamma=1$, community-based variants of the Barabási–Albert model have also been investigated. In particular,~\cite{bhamidi2022community} studies the case of our \model \ 1 with $q_A=q_B=q$ and $\gamma\in\{0, 1\}$. Since individual reproduction rates depend on the overall population composition, the resulting process is not a classical branching process. By combining continuous-time embedding techniques with stochastic analysis, the authors establish the convergence of certain macroscopic graph quantities and use these results to derive further asymptotic properties. In addition, Jordan considers the case of finitely many communities in~\cite{jordan2013geometric}, and even extends the framework to a compact metric space of communities, still under linear preferential attachment. More generally, he studies an extension of the model in which each new node arrives with a fixed number $m \in \N^*$ of edges. His approach consists in first establishing a convergence result for the analogues of the ratios $\frac{D_n(A)}{n}$ and $\frac{D_n(B)}{n}$ in our model, using an initial stochastic approximation result, as well as on the existence of a suitably chosen Lyapunov function. This convergence result is then used to derive a new recurrence relation for the proportion of vertices of given degrees and types in the sequence of graphs $(G_n)_{n \in \N}$. Finally, a second stochastic approximation result is applied to establish the desired convergence. 

The difficulty in extending the approach of Jordan lies in the need to establish beforehand the convergence of the ratios $\frac{D_n(A)}{n}$ and $\frac{D_n(B)}{n}$. It is not clear to us whether a Lyapunov function satisfying the definition given in~\cite{jordan2013geometric} can be found in our setting. This is why we found it necessary to develop a new approach to address the sublinear multitype case, which consists in establishing a convergence result, as $n$ tends to infinity, for the counting measure $m_n$ describing the proportion of vertices associated with each pair $(k, \Lambda)\in \N^* \times \{A, B\}$. 

Furthermore, the paper~\cite{antunes2013attribute} studies a class of attribute-based preferential attachment models, in which new vertices join the network sequentially and connect to existing vertices according to their degrees and attributes, which play a role analogous to community labels in community-based network models. The authors establish asymptotic results on the local structure of these networks, the tail behavior of degree and PageRank score distributions, and the performance of sampling methods in capturing different attribute groups. Their analysis combines local weak convergence, continuous-time multitype branching processes, and stochastic approximation techniques (drawing in particular on the work of Jordan in~\cite{jordan2013geometric}) to characterize the limiting network structure and assess sampling biases. However, most of this work focuses on the case where the weight function is linear.

\begin{remark}
For the study of our two-community model, one could similarly introduce the following continuous-time multitype branching process. Without loss of generality, assume that $G_0$ consists of a single node, $\emptyset$. This initial node is assigned type $A$ with probability $p$ and type $B$ with probability $1-p$. The initial node gives birth to vertices of types $A$ and $B$ at time-dependent rates to be specified just below. Similarly, each newly born node, in turn, gives birth to new vertices of types $A$ and $B$ at the following time-dependent rates. Let $G_t$ denote the graph obtained at time $t$ and let $V(G_t)$ be its node set. In \model \ $1$, a node $v \in V(G_t)$ of type $\Lambda := \sigma(v)$ gives birth to a node of type $\Lambda'$ at rate 
\[
r^{(1)}_{v \to \Lambda'}(t) = \cfrac{p_{\Lambda'}\times[q_{\Lambda'} \ind_{\Lambda \neq \Lambda'} + (1-q_{\Lambda'})\ind_{\Lambda = \Lambda'}] \times \deg_t(v)^\gamma}{D_t(\Lambda)}\, ,
\]
where $\deg_t(v)$ denotes the degree of $v$ in $G_t$ and $D_t(\Lambda)$ the following continuous-time analogue of $D_n(\Lambda)$ 
\[
D_t(\Lambda) := \displaystyle\sum_{\substack{v \in V(G_t) \\ \sigma(v) = \Lambda}} \deg_t(v)^\gamma \, .
\]
In \model \ $2$, we have 
\[
r^{(2)}_{v \to \Lambda'}(t) = \cfrac{p_{\Lambda'} \times [q_{\Lambda'} \ind_{\Lambda \neq \Lambda'} + (1-q_{\Lambda'})\ind_{\Lambda = \Lambda'}] \times \deg_t(v)^\gamma}{(1-q_{\Lambda'})D_t(\Lambda') + q_{\Lambda'}D_t(\overline{\Lambda'})} \, .
\]
We then define the sequence of stopping times $(T_n)_{n\geq 0}$, corresponding to the times at which a new node appears in the randomly growing tree
\[
T_n := \inf\{t \geq 0 : \vert V(G_t) \vert = n +1 \} \, .
\]
With this notation, the sequence of growing graphs $(G_n^{\mathrm{d}})_{n \geq 0}$ generated by the discrete-time model defined in Definition~\ref{def:model-1} is related to the continuous-time model $(G_t)_{t \geq 0}$ introduced above through the following equality in distribution:
\[
G_n^{\mathrm{d}} \overset{(\mathrm{d})}{=} G_{T_n} \, .
\]
Yet a key difference arises at this point: the point processes associated with each node of the continuous-time growing graph, which describe the birth times and types of its offspring, depend both on the type of the parent node and on the structure of the entire graph. To generalize the approach used in~\cite{rudas2007random}, one could rely on multi-type extensions of the convergence results employed therein. In this regard, a paper by Jagers and Nerman seems particularly relevant (see~\cite{jagers1996asymptotic}). However, to apply the result of this paper, as well as to extend the approach of~\cite{rudas2007random} to the two-community setting, one would first need to establish a convergence result for the ratios $\frac{D_n(\Lambda)}{D_n(\Lambda)+D_n(\overline{\Lambda})}$ or $\frac{D_n(\Lambda)}{n}$, for $\Lambda \in \{A, B\}$. The continuous-time embedding of our model cannot be interpreted as a branching process in the usual sense (as defined, for instance, in~\cite{jagers1996asymptotic}), since the individual reproduction rates at any given time $t$ depend on the composition of the entire population at that time. This obstacle is all the more interesting in that it also arises when attempting to generalize the approach used by Jordan in~\cite{jordan2013geometric} where $\gamma =1$. In addition, it is worth noting that, Theorem~\ref{th:deg-distrib-2} is the key step in justifying the approximation of the continuous-time embedding introduced above by a branching process whose rates at time $t$ no longer depend on the entire graph structure at that time. We emphasize that the rates of this approximating branching process are obtained by multiplying the original rates by $D_t(\Lambda) + D_t(\overline{\Lambda})$, which corresponds to a time change of the continuous-time process $(G_t)$, and then replacing the ratios $\frac{D_t(\Lambda)}{D_t(\Lambda)+D_t(\overline{\Lambda})}$ and $\frac{D_t(\overline{\Lambda})}{D_t(\Lambda)+D_t(\overline{\Lambda})}$ appearing in the resulting rates with their respective limits.
\end{remark}

\subsection{Proof strategy and outline of the paper}

\subsubsection{Measure-valued P\'olya processes}

The process $(m_n)_{n \in \N}$ can be described using the formalism of \emph{measure-valued P\'olya processes}, as studied in~\cite{mailler2017measure} and~\cite{mailler2020stochastic}. These processes provide a generalization of P\'olya urn processes by allowing the number of colors to be infinite, with values in a Polish space $E$. In our setting, starting from a measure $m_0$ encoding the number of nodes in $G_0$ with a given degree and a given type, we obtain $m_{n+1}$ from $m_n$ in the following way. At time $n+1$, the new node $v_{n+1}$ gets a type $\Lambda_{n+1}$ which is equal to $A$ with probability $p$, or to $B$ with probability $1-p$. We let $B_{n+1}$ be the Bernoulli random variable given by $B_{n+1}=\ind_{\{\Lambda_{n+1}=A\}}$. Letting $X_{n+1} \in E$ the pair corresponding to the degree and type of the node to which $v_{n+1}$ attaches to, we have, in \model\ $1$, for all $x\in E$, and $\Lambda\in\{A, B\}$,
\[
\P(X_{n+1} = x \mid m_n, \Lambda_{n+1}=\Lambda) = (1-q_\Lambda)\frac{m_nD^\Lambda(x)}{m_n D^\Lambda(E)}+q_\Lambda \frac{m_nD^{\overline{\Lambda}}(x)}{m_n D^{\overline{\Lambda}}(E)}\, ,
\]
where $D^A$ is the infinite diagonal matrix with diagonal entries: $1, 0, 2^\gamma, 0, 3^\gamma, 0...$, and $D^B$ is the infinite diagonal matrix with diagonal entries: $0, 1, 0, 2^\gamma, 0, 3^\gamma...$. We specify that if $\eta$ is a measure on $E$ and $M$ is an infinite matrix in $\mathcal{M}_{E\times E}(\mathbb{R})$, then for $x\in E$, we have $\eta M(x)= \sum_{y \in E} \eta(y) M(y, x)$, and $\eta M(E)=\sum_{x\in E}\eta M(x)$.

In \model\ $2$, we have for all $x\in E$ and $\Lambda\in\{A, B\}$,
\[ 
\P(X_{n+1} = x \mid m_n, \Lambda_{n+1}=\Lambda) = \cfrac{m_nP^\Lambda(x)}{m_nP^\Lambda (E)} \, ,
\]
where $P^\Lambda=(1-q_\Lambda)D^\Lambda +q_\Lambda D^{\overline{\Lambda}}$. 

Put differently, in \model\ $2$, setting $P^{(n+1)}=B_{n+1}P^A+(1-B_{n+1})P^B$, we have
\[
\P(X_{n+1}=x \mid m_n, \Lambda_{n+1}) =\frac{m_n P^{(n+1)}(x)}{m_n P^{(n+1)}(E)}\, .
\]

In both models, the replacement rule is then: $m_{n+1} = m_n + R^{(n+1)}_{X_{n+1}}$, with
\[
R^{(n+1)}_{(k, \Lambda)} = \delta_{(1, \Lambda_{n+1})} - \delta_{(k, \Lambda)} + \delta_{(k+1, \Lambda)}\, .
\]
The matrix $R^{(n+1)}$, whose rows are $(R^{(n+1)}_{x})_{x\in E}$, is called the replacement matrix. As in~\cite{mailler2020stochastic}, the replacement matrices are random: $R^{(n+1)}$ depends on the random type $\Lambda_{n+1}$.

\begin{remark}
While we use the formalism of measure-valued P\'olya processes presented in~\cite{mailler2020stochastic}, it is important to note that none of the two models falls into the framework of this paper, in which the weight matrices are all equal to the same deterministic matrix $P$. \model\ $2$ could be seen as a generalization where the weight matrices $P^{(n)}$ are random (with the same source of randomness as the one of the replacement matrices $R^{(n)}$).
\end{remark}

\subsubsection{Stochastic approximation}\label{sec:outline-stochastic-approx}

Rather than studying the process $(m_n)$ directly, it is convenient to first focus on the process $(\eta_n)_{n\geq 0}$, defined by $\eta_0=0$ and for $n\geq 1$, $\eta_n=\frac{1}{n}\sum_{i=1}^n \delta_{X_i}$. By definition, we have
\[
\eta_{n+1}=\left(1-\frac{1}{n+1}\right)\eta_n+\frac{1}{n+1}\delta_{X_{n+1}}\, .
\]
From now on, we focus on \model \ $1$. Observe that
\[
\E\left[\delta_{X_{n+1}}\given \cF_n\right]= r_A \frac{m_n D^A}{m_n D^A(E)}+ r_B\frac{m_n D^B}{m_n D^B(E)} \, ,
\]
where $r_\Lambda=p_\Lambda(1-q_\Lambda)+p_{\overline{\Lambda}}q_{\overline{\Lambda}}$, and that
\begin{equation}\label{eq:m_n-eta}
m_n=m_0+\sum_{i=1}^n R_{X_i}^{(i)}=m_0+\sum_{i=1}^n \delta_{X_i}R^{(i)}\, ,
\end{equation}
where 
\[
R^{(i)}=\begin{pmatrix}
B_i & 1-B_i & 0 & \dots & 0 & \dots \\
\vdots & \vdots & \vdots & & \vdots & \\
B_i & 1-B_i & 0 & \dots & 0 & \dots \\
\vdots & \vdots & \vdots & & \vdots &
\end{pmatrix}
+
\begin{pmatrix}
-1& 0 & 1 & 0 &  & (0) & \\
 & \ddots & \ddots & \ddots & \ddots & &\\
 &  & -1 & 0 & 1 & 0 & \\
 & (0) & & \ddots & \ddots & \ddots & \ddots
\end{pmatrix}\, .
\]
If in~\eqref{eq:m_n-eta} each matrix $R^{(i)}$ could be approximated by its expectation $R$, where $B_i$ is replaced by $p$, we would have
\[
\frac{m_n}{n}\approx \eta_n R\, .
\]
Defining $Q^A=RD^A$ and $Q^B=RD^B$, we would get 
\[
\E\left[\delta_{X_{n+1}}\given \cF_n\right]\approx r_A \frac{\eta_n Q^A}{\eta_n Q^A(E)}+r_B\frac{\eta_n Q^B}{\eta_n Q^B(E)}\, .
\]
This heuristic leads us to write
\[
\eta_{n+1}=\eta_n +\frac{1}{n+1}\left(r_A \frac{\eta_n Q^A}{\eta_n Q^A(E)}+r_B\frac{\eta_n Q^B}{\eta_n Q^B(E)} -\eta_n +\varepsilon_{n+1} \right) \, ,
\]
where 
\[
\varepsilon_{n+1}= \delta_{X_{n+1}}-r_A \frac{\eta_n Q^A}{\eta_n Q^A(E)} -r_B\frac{\eta_n Q^B}{\eta_n Q^B(E)}\, .
\]
In Section~\ref{sec:stochastic-approximation}, we show, using stochastic approximation techniques developed in~\cite{benaim1999dynamics}, that for all $\gamma\in [0,1)$, an affine interpolation of $(\eta_n)_{n\geq 0}$ with appropriate step sizes is an \emph{asymptotic pseudo-trajectory} of the semiflow on $\cP(E)$ (the set of probability distributions on $E$) induced by
\[
\partial_t \eta_t =F(\eta_t)\, ,
\]
for some explicit function $F$ from $\cP(E)$ to $\cM_0(E)$, the set of signed measures on $E$ with zero sum (Proposition~\ref{prop:semiflow1}). 

In Section~\ref{sec:fixed-point-unique}, we show that there exists a unique fixed point to the continuous-time dynamics induced by $F$, \emph{i.e.}\ a unique probability distribution $\eta^\star\in\cP(E)$ satisfying $F(\eta^\star)=0$ (Proposition~\ref{prop:stationnarysolution}). The next task is to show that $\eta_t$ converges to $\eta^\star$ as $t$ tends to infinity.

\subsubsection{Non-linear Markov processes and coupling}

Since the dynamics induced by $F$ preserves probability distributions, we would like to view it as distribution of some stochastic process on $E$. However, no such process is immediately given from $F$ itself. Indeed, we have $F(\eta)=\eta A_\eta$ but the matrix $A_\eta$ is not an infinitesimal generator because the sum on each row is not equal to zero. Still, one may design a small perturbation $L_\eta$ of $A_\eta$, which is such that $F(\eta)=\eta L_\eta$, which has negative entries on the diagonal and which has zero sum on each row. In the two-community case, to ensure that all the off-diagonal entries of $L_\eta$ are positive, we could not do without a condition on $\gamma$ (see Assumption~\ref{assume:gamma}), but such a condition is not needed in the one-community setting where the corresponding matrix $L_\eta$ corresponds to an infinitesimal generator for all $\gamma\in [0, 1]$. In Section~\ref{sec:one-community}, we first focus on this last setting, which allows us to introduce in a simpler setting the arguments that will be used for two communities. Having described our dynamics as the one of a non-linear Markov process on $\N^*$ (the non-linearity coming from the fact that the generator $L_\eta$ depends on the measure $\eta$ itself), we use coupling methods inspired by~\cite{monmarche2023elementary} to establish that, for any starting distribution, the equation $\partial_t\eta_t=\eta_t L_{\eta_t}$ has a unique solution on $\R_+$ (Proposition~\ref{prop:existence-uniqueness-one-community}), and that, whatever the starting distribution, this solution converges to $\eta^\star$ in total variation (Proposition~\ref{prop:coupling-one-community}). This last result holds for all $\gamma\in [0, 1]$ satisfying inequality~\eqref{eq:inequality}. We strongly believe that this inequality holds for all $\gamma\in [0, 1)$ (see Figure~\ref{fig:function-inequality}), but we only have a simple analytic proof for $\gamma \in \left[0, \tfrac{\ln 2}{\ln 3}\right]$. In Section~\ref{sec:convergence-two-communities}, we come back to the two-community setting. We first use a fixed point theorem to establish existence and uniqueness of a solution to $\partial_t\eta_t=\eta_t L_{\eta_t}$ (Proposition~\ref{prop:semiflow}). Then we adapt the coupling arguments that will be introduced in Section~\ref{sec:one-community} to the case of two communities and show that, for any starting distribution, the unique solution to $\partial_t\eta_t=\eta_t L_{\eta_t}$ converges to $\eta^\star$. Once this is done, the proof of Theorem~\ref{th:deg-distrib-2} (and of Theorem~\ref{th:deg-distrib-1}) easily follows. Section~\ref{sec:model2} highlights the few changes that have to be made to adapt the proof to \model\ $2$.

\section{Stochastic approximation}\label{sec:stochastic-approximation}
Recall from Section~\ref{sec:outline-stochastic-approx} that the approximation $\tfrac{m_n}{n}\approx \eta_n R$ led us to write
\begin{equation}\label{eq:dynamic-1}
\eta_{n+1}=\eta_n +\frac{1}{n+1}\left(r_A \frac{\eta_n Q^A}{\eta_n Q^A(E)}+r_B\frac{\eta_n Q^B}{\eta_n Q^B(E)} -\eta_n +\varepsilon_{n+1} \right) \, ,
\end{equation}
where 
\[
\varepsilon_{n+1}= \delta_{X_{n+1}}-r_A \frac{\eta_n Q^A}{\eta_n Q^A(E)} -r_B\frac{\eta_n Q^B}{\eta_n Q^B(E)}\, .
\]
The next lemma, which quantifies the approximation $\tfrac{m_n}{n}\approx \eta_n R$, will be used repeatedly.

\begin{lemma}\label{lem:m_n-eta_n}
For all $n\geq 1$, we have
\[
\frac{m_n}{n}-\eta_nR=\frac{m_0}{n}+\frac{1}{n}\sum_{i=1}^n(B_i-p, p-B_i, 0, 0, \dots)\, .
\]
\end{lemma}
\begin{proof}
Using~\eqref{eq:m_n-eta}, we have
\begin{align*}
\frac{m_n}{n}-\eta_nR & = \frac{m_0}{n}+\frac{1}{n} \sum_{i=1}^n \delta_{X_i}(R^{(i)}-R) \\
    &= \frac{m_0}{n}+\frac{1}{n} \sum_{i=1}^n \delta_{X_i}\begin{pmatrix}
B_j-p & p-B_j & 0 & \dots & 0 & \dots \\
\vdots & \vdots & \vdots & & \vdots & \\
B_j-p & p-B_j & 0 & \dots & 0 & \dots \\
\vdots & \vdots & \vdots & & \vdots &
\end{pmatrix}  \\
&= \frac{m_0}{n}+\frac{1}{n} \sum_{i=1}^n \left(B_i-p, p-B_i\dots, 0, \dots\right)\, ,
\end{align*}
which concludes the proof.
\end{proof}
It will be convenient to slightly change the step sizes and to rewrite Equation~\eqref{eq:dynamic-1} in the form
\begin{align*}
    \eta_{n+1} &= \eta_n +\gamma_{n+1}\left\{F(\eta_n)+\varepsilon'_{n+1}\right\}\, ,
\end{align*}
where
\begin{align*}
\gamma_{n+1} &= \cfrac{1}{n+1} \left(\cfrac{r_A}{\eta_n Q^A(E)} + \cfrac{r_B}{\eta_nQ^B(E)} \right)\, ,\\
\varepsilon'_{n+1}&=\cfrac{\varepsilon_{n+1}}{\cfrac{r_A}{\eta_n Q^A(E)} + \cfrac{r_B}{\eta_nQ^B(E)}}\, ,
\end{align*}
and where
\begin{align}\label{eq:def-F}
F(\eta)&=s(\eta)\left(\eta Q^A-\eta Q^A(E)\eta\right)+(1-s(\eta))\left(\eta Q^B-\eta Q^B(E)\eta\right)\, ,
\end{align}
with
\begin{equation}\label{eq::def-s}
s(\eta) = \cfrac{r_A\eta Q^B(E)}{r_B\eta Q^A(E) + r_A\eta Q^B(E)}\, .
\end{equation}
We introduce the following notations for later use: 
\[
s^A(\eta) = s(\eta) \quad \text{and} \quad s^B(\eta) = 1 - s(\eta).
\]
Note that $\gamma_{n+1}$ is of the order $\tfrac{1}{n+1}$. Indeed, for all $\eta\in \cP(E)$ and for $\Lambda\in\{A, B\}$, we have
\[
p_\Lambda\leq \eta Q^\Lambda(E) = p_\Lambda+\sum_{k\geq 1}\eta(k, \Lambda)\left((k+1)^\gamma-k^\gamma\right)\leq p_\Lambda+2^\gamma-1\, .
\]
We then consider the increasing sequence $(\tau_n)_{n \geq 0}$ defined by $\tau_0=0$ and for all $n\geq 1$,
\[
\tau_n = \gamma_1 + \gamma_2 + ... + \gamma_n\, ,
\]
and we introduce the linearized version $(X(t))_{t\geq 0}$ of $(\eta_n)_{n \in \N}$ defined, for all $n \geq 1$ and all $t \in [\tau_n, \tau_{n+1}]$, by
\[
X(t) = \eta_n + \cfrac{t-\tau_n}{\tau_{n+1}-\tau_n}(\eta_{n+1}-\eta_n)\, .
\]

Now, let $\Phi$ be the semiflow induced by the dynamical system
\begin{align}\label{eq:dynamicsystem-1}
    \partial_t \eta_t &= F(\eta_t)\,,
\end{align}
with $F$ defined in~\eqref{eq:def-F}, that is, the semiflow given by
\begin{eqnarray*}
\Phi\colon \R_+\times\cP(E) &\to & \cP(E)\\
 (t, \eta) &\mapsto & \Phi(t, \eta)=\Phi_t(\eta)=\eta_t\, ,
\end{eqnarray*}
where $\eta_t$ is the state at time $t$ of the system given by $\partial_t\eta_t=F(\eta_t)$, with initial condition $\eta$. The fact that this semiflow is well-defined will be established in Proposition~\ref{prop:semiflow}. The rest of this section is devoted to the proof of the following result.

\begin{prop}\label{prop:semiflow1}
For all $\gamma\in [0, 1)$, the interpolated process $(X(t))_{t\geq 0}$ is an asymptotic pseudo-trajectory of the semiflow $\Phi$ induced by $F$, that is, for all $T>0$, we have
\[
\lim_{t\to +\infty} \sup_{0\leq h\leq T}\| X(t+h)-\Phi_h(X(t))\|_1=0\, .
\]
\end{prop}

To prove Proposition~\ref{prop:semiflow1}, we apply \cite{benaim1999dynamics}[Theorem 3.2] which reduces the problem to showing three points:
\begin{enumerate}
    \item uniform continuity: the function $X:\R_+ \to \cP(E)$ is uniformly continuous;
    \item compact closure: the image of $X:\R_+ \to \cP(E)$ has compact closure in $\cP(E)$;
    \item every limit point of the translation flow\footnote{The definition of the translation flow will be given in Section~\ref{sec:limit-translation}.} applied to $X$ belongs to the set 
    \[
    \{t\mapsto \Phi_t(\eta),\,\eta\in\cP(E)\}\, .
    \]
\end{enumerate}
Point~$(1)$ is established in Section~\ref{sec:unif_continue}, Point~$(2)$ in Section~\ref{sec:compact-closure}, Point~$(3)$ in Section~\ref{sec:limit-translation}.

\subsection{Uniform continuity}\label{sec:unif_continue}

\begin{lemma}\label{lem:unif_continue}
The function $X:\R_+ \to \cP(E)$ is uniformly continuous.
\end{lemma}
\begin{proof}
For $t\geq 0$, let us define
\[
m(t):=\sup\{k\geq 0,\, t\geq \tau_k\}\, .
\]
Using a slight reformulation of~\cite{benaim1999dynamics}[Equation (9)], we have, for all $s, t\geq 0$,
\[
X(t+s)-X(t)=\int_t^{t+s} \frac{1}{\gamma_{m(u)+1}(m(u)+1)}\left(\delta_{X_{m(u)+1}}-\eta_{m(u)}\right) d s\, .
\]
Using that $\gamma_{n+1}(n+1)$ is uniformly bounded away from $0$, and that both $\delta_{X_{n+1}}$ and $\eta_n$ are probability distributions, we obtain that there exists $c>0$ such that for all $s, t\geq 0$,
\[
\|X(t+s)-X(t)\|_1\leq cs\, .
\]
In particular, the function $X$ is uniformly continuous.
\end{proof}

\subsection{Compact closure}\label{sec:compact-closure}

To prove that $\{X(t),\, t\in\R_+\}$ has compact closure, let us show the following.

\begin{lemma}\label{lem:compact-closure}
For all $\gamma\in [0,1[$, there exists $\delta>0$ small enough, and a random value $C\in \R_+$, such that, defining $W: E\to \N^*$ by
\[
\forall (k, \Lambda)\in E\, , \; W(k, \Lambda)=k^\delta\, .
\]
we have
\[
\forall n\geq 1\, ,\; \langle \eta_n, W\rangle :=\sum_{x\in E}\eta_n(x)W(x) \leq C\, .
\]
\end{lemma}

\begin{proof}
By~\eqref{eq:dynamic-1}, letting $u_n=\langle \eta_n, W\rangle$, we have, for all $n\geq 1$,
\[
u_{n+1} =u_n+\frac{1}{n+1}\left(r_A\left\langle \frac{\eta_n Q^A}{\eta_n Q^A(E)}, W\right\rangle +r_B\left\langle \frac{\eta_n Q^B}{\eta_n Q^B(E)}, W\right\rangle -u_n +\langle \varepsilon_{n+1}, W\rangle\right)\, .
\]
On the one hand, if $\delta>0$ is chosen small enough so that $\gamma+\delta<1$, then for $\Lambda\in \{A, B\}$, we have
\begin{align*}
\left\langle \frac{\eta_n Q^\Lambda}{\eta_n Q^\Lambda(E)}, W\right\rangle & =\frac{1}{\eta_n Q^\Lambda(E)}\langle \eta_n, Q^\Lambda W\rangle\\
&= \frac{1}{\eta_n Q^\Lambda(E)}\left(p_\Lambda+\sum_{k\geq 1}\eta_n(k, \Lambda)\left((k+1)^{\gamma+\delta}-k^{\gamma+\delta}\right)\right)\\
&\leq \frac{p_\Lambda+1}{p_\Lambda}\, .
\end{align*}
Hence there exists a real-valued constant $a>0$ such that
\[
r_A\left\langle \frac{\eta_n Q^A}{\eta_n Q^A(E)}, W\right\rangle +r_B\left\langle \frac{\eta_n Q^B}{\eta_n Q^B(E)}, W\right\rangle\leq a\, .
\] 
We thus have
\[
u_{n+1} \leq u_n+\frac{1}{n+1}\left(a -u_n +\langle \varepsilon_{n+1}, W\rangle\right)\, .
\]
If we could show that 
\[
L:=\sup_{1\leq n\leq N}\left|\sum_{i=n}^N \frac{1}{i+1}\langle \varepsilon_{i+1}, W\rangle\right| 
\]
is finite almost surely, then we would get that for all $n\geq 1$, $u_n \leq 2(L+ c)$, with $c=a\vee u_1$. Indeed, this is immediate for $n=1$ and, assuming that it is true for some $n\geq 1$:
\begin{itemize}
    \item if $u_n\leq c$, then $u_n+1\leq c+a+L$ and the result holds for $n+1$;
    \item if $u_n> c$, then define $n_0$ by
    \[
    n_0=\sup\left\{k\in\{1,\dots, n\}:\, u_k>c,\, u_{k-1}\leq c\right\}\, ,
    \]
    which is well defined because $u_1\leq c$. By the above, we have $u_{n_0}\leq c+a+L$, and, since $u_{n_0}, \dots, u_n>c$, we have
    \[
    u_n\leq u_{n_0}+\sum_{i=n_0}^n \frac{1}{i+1}\langle \varepsilon_{i+1}, W\rangle\leq c+a+2L\, ,
    \]
    and the result holds for $n+1$.
\end{itemize}
To prove that $L$ is finite almost surely, let us show that it is bounded in $L^1$. We have
\begin{align*}
\E[L] &\leq \underbrace{\E\left[\sup_{1\leq n\leq N}\left| \sum_{i=n}^N \frac{1}{i+1}\langle \delta_{X_{i+1}} -\E[\delta_{X_{i+1}}|\cF_i], W\rangle\right|\right]}_{(\texttt{I})}\\
&\hspace{1cm}+ \underbrace{\sum_{\Lambda\in\{A, B\}}r_\Lambda\sum_{i=1}^{+\infty}\frac{1}{i+1}\E\left[\left|\left\langle \frac{m_iD^\Lambda}{m_iD^\Lambda(E)}-\frac{\eta_iQ^\Lambda}{\eta_iQ^\Lambda(E)}, W\right\rangle\right|\right]}_{(\texttt{II})}\, .
\end{align*}
For term $(\texttt{I})$, using the Jensen Inequality and the Burkholder--Davis--Gundy Inequality, we see that there exists $C_2>0$ such that
\begin{equation}\label{eq:bound-L-I}
    (\texttt{I}) \leq \sqrt{C_2\sum_{i=1}^{+\infty} \frac{1}{(i+1)^2}\E\left[\left(W(X_{i+1}) -\E[W(X_{i+1})|\cF_i]\right)^2\right]} \lesssim \sqrt{\sum_{i=1}^{+\infty} \frac{1}{(i+1)^{2-2\delta}}}\, ,
\end{equation}
where we used that $W(X_{i+1})$ is deterministically smaller, up to constants\footnote{Here and in the rest of the paper, for two functions $f$ and $g$ from $\N^*$ to $\R_+$, we write $f(n)\lesssim g(n)$ to say that there exists a constant $C\in (0,+\infty)$ such that for all $n\in\N^*$, we have $f(n)\leq Cg(n)$.}, than $i^\delta$. Hence, taking $\delta<1/2$ ensures that $(\texttt{I})$ is bounded. Moving to term $(\texttt{II})$, we set $\tilde{m}_i=\tfrac{m_i}{i}$. Clearly, $\tfrac{m_i D^\Lambda}{m_i D^\Lambda(E)}=\tfrac{\tilde{m}_i D^\Lambda}{\tilde{m}_i D^\Lambda(E)}$, and by the triangle inequality, we have, for all $i\geq 1$ and $\Lambda\in\{A, B\}$, 
\begin{align*}
\E\left[\left|\left\langle \frac{\tilde{m}_iD^\Lambda}{\tilde{m}_iD^\Lambda(E)}-\frac{\eta_iQ^\Lambda}{\eta_iQ^\Lambda(E)}, W\right\rangle\right|\right] &\leq
\underbrace{\E\left[\left|\left\langle \frac{\tilde{m}_iD^\Lambda-\eta_iQ^\Lambda}{\eta_iQ^\Lambda(E)}, W\right\rangle\right|\right]}_{(\texttt{II.a})}\\
&\hspace{2cm}+\underbrace{\E\left[\left|\left\langle \frac{\tilde{m}_iD^\Lambda}{\tilde{m}_iD^\Lambda(E)}-\frac{\tilde{m}_i D^\Lambda}{\eta_iQ^\Lambda(E)}, W\right\rangle\right|\right]}_{(\texttt{II.b})}\, .
\end{align*}
For term $(\texttt{II.a})$, using the facts that $\eta_iQ^\Lambda(E)$ is lower bounded by $p_\Lambda$, that $\tilde{m}_iD^\Lambda-\eta_iQ^\Lambda=(\tilde{m}_i-\eta_i R)D^\Lambda$, and applying Lemma~\ref{lem:m_n-eta_n}, we get
\begin{equation}\label{eq:bound-L-IIa}
(\texttt{II.a}) \leq \frac{W(1, \Lambda)}{p_\Lambda}\left(\frac{m_0(1, \Lambda)}{i}+\frac{1}{i}\E\left[\left|\sum_{j=1}^i (B_j-p)\right|\right]\right)\lesssim\frac{1}{\sqrt{i}}\, ,
\end{equation}
where we used that $\E\left[\left|\sum_{j=1}^i (B_j-p)\right|\right]\leq \sqrt{p(1-p)i}$. For term $(\texttt{II.b})$, we write
\[
(\texttt{II.b})=\E\left[\left|1-\frac{\tilde{m}_iD^\Lambda(E)}{\eta_i Q^\Lambda(E)}\right|\left\langle \frac{\tilde{m}_iD^\Lambda}{\tilde{m}_iD^\Lambda(E)}, W\right\rangle\right]\lesssim \E\left[\left|\eta_i Q^\Lambda(E)-\tilde{m}_iD^\Lambda(E)\right|\right]i^\delta \, ,
\]
where we used that $\eta_i Q^\Lambda(E)$ is larger than $p_\Lambda$ and that the probability measure $\tfrac{\tilde{m}_iD^\Lambda}{\tilde{m}_iD^\Lambda(E)}$ is supported on degrees at most $i+d_{\textrm{max}}$, where $d_{\textrm{max}}$ is the maximum degree in $G_0$. Hence, another application of Lemma~\ref{lem:m_n-eta_n} yields 
\begin{equation}\label{eq:bound-L-IIb}
(\texttt{II.b})\lesssim \frac{1}{i^{1/2-\delta}}\, .
\end{equation}
Combining~\eqref{eq:bound-L-I}, \eqref{eq:bound-L-IIa} and \eqref{eq:bound-L-IIb}, we obtain
\[
\E[L]\lesssim \sqrt{\sum_{i=1}^{+\infty} \frac{1}{i^{2-2\delta}}}+\sum_{i=1}^{+\infty}\frac{1}{i^{3/2}}+\frac{1}{i^{3/2-\delta}}\, ,
\]
which is finite for $\delta<1/2$, and thus completes the proof.
\end{proof}

\subsection{Limit points of the translation flow are solutions}\label{sec:limit-translation}

The translation flow $\Theta: C^0(\R_+, \cP(E))\times\R_+\to C^0(\R_+, \cP(E))$ is defined by
\[
\Theta^t(X)(s)\mapsto X(t+s)\, .
\]

\begin{prop}\label{prop:limit-translation}
For all $\gamma\in [0,1)$, every limit point of the sequence $\{\Theta^t(X)\}_{t\geq 0}$ belongs to the set $\{t\mapsto \Phi_t(\eta),\,\eta\in\cP(E)\}$.
\end{prop}

Before proving Proposition~\ref{prop:limit-translation}, let us introduce the same notation as in~\cite{benaim1999dynamics}. For $u\in [\tau_n,\tau_{n+1}[$, we define
\[
\overline{X}(u)=\eta_n \quad \text{ and }\quad \overline{\varepsilon}(u)=\varepsilon'_{n+1}\, .
\]
\begin{proof}[Proof of Proposition~\ref{prop:limit-translation}]
We proceed as in the proof of~\cite{benaim1999dynamics}[Proposition 4.1] and decompose $\Theta^t(X)$ as
\[
\Theta^t(X)=L_F(\Theta^t(X))+A_t+B_t\, ,
\]
where $L_F: C^0(\R_+, \R^E)\to C^0(\R_+, \R^E)$ is given by
\[
L_F(X)(s)=X(0)+\int_0^s F(X(u))\d u\, ,
\]
and where
\begin{eqnarray*}
A_t(s) &=&  \int_t^{t+s} \{F(\overline{X}(u))-F(X(u))\}\d u \, , \\
B_t(s) &=& \int_t^{t+s} \overline{\varepsilon}(u)\d u\, .
\end{eqnarray*}
We will show that\footnote{This is more than enough: it would suffice to have it for the supremum over $0\leq s\leq T$ for any fixed $T>0$.}
\begin{equation}\label{eq:A_t-B_t}
    \sup_{s\geq 0}\|A_t(s)\|_1 \, \underset{t\to +\infty}\longrightarrow\, 0 \quad\text{ and} \quad \sup_{s\geq 0}\|B_t(s)\|_1\, \underset{t\to +\infty}\longrightarrow\, 0\, .
\end{equation}
Let us first control the term $A_t$. Note that for some measure $\eta\in\cP(E)$, we can write
\[
F(\eta)=s(\eta)\eta Q^A+(1-s(\eta))\eta Q^B -\lambda(\eta)\eta\, ,
\]
where
\[
\lambda(\eta)=\frac{\eta Q^A(E)\eta Q^B(E)}{r_A\eta Q^B(E)+r_B\eta Q^A(E)}\, \cdot 
\]
With this writing at hand, repeated applications of the triangle inequality gives
\begin{align*}
    \| F(\eta)-F(\nu)\|_1&\leq s(\eta)\|(\eta-\nu)Q^A\|_1 +(1-s(\eta))\|(\eta-\nu)Q^B\|_1 \\
    &\hspace{0.5cm}+|s(\eta)-s(\nu)|\left(\|\nu Q^A\|_1+\|\nu Q^B\|_1 \right)+ \lambda(\eta)\|\eta-\nu\|_1 +|\lambda(\eta)-\lambda(\nu)|\, .
\end{align*}
We claim that
\begin{equation}\label{eq:ineq-s(eta)-tv}
    |s(\eta)-s(\nu)|\leq \frac{2^\gamma(2^\gamma-1)}{8p_Ap_B}\|\eta-\nu\|_{1}\, ,
\end{equation}
and that
\begin{equation}\label{eq:ineq-c(eta)-tv}
    |\lambda(\eta)-\lambda(\nu)|\leq \frac{1}{2}\left(1+\frac{1}{2p_Ap_B}\right)\|\eta-\nu\|_{1}\, .
\end{equation}
The proof of those two inequalities can be found in Section~\ref{sec:proof-ineq-s(eta)-c(eta)-tv} (for~\eqref{eq:ineq-s(eta)-tv}, the precise constant factor is not useful here but it will be used later in the paper). Omitting the precise constant factors, we have
\[
 \| F(\eta)-F(\nu)\|_1\lesssim \|\eta-\nu\|_1 +\max_{\Lambda\in \{A, B\}}\left\{ \|(\eta-\nu)Q^\Lambda\|_1 +\|\eta-\nu\|_1\|\nu Q^\Lambda\|_1\right\}\, .
\]
Choosing $\eta=X(u)$ and $\nu=\overline{X}(u)$ with $u\in [\tau_n,\tau_{n+1}[$ and using that $\eta_{n+1}-\eta_n=\tfrac{1}{n+1}(\delta_{X_{n+1}}-\eta_n)$, this yields
\begin{align*}
\| F(X(u))-F(\overline{X}(u))\|_1
&\lesssim \|\eta_{n+1}-\eta_n\|_1 \\
&\qquad\qquad + \max_{\Lambda\in \{A, B\}}
\left\{\|(\eta_{n+1}-\eta_n)Q^\Lambda\|_1 +\|\eta_{n+1}-\eta_n\|_1 \|\eta_n Q^\Lambda\|_1 \right\} \\
&\lesssim \frac{1}{n+1} \left( 1 + \max_{\Lambda\in \{A, B\}} \left\{\|\delta_{X_{n+1}}Q^\Lambda \|_1 +\|\eta_n Q^\Lambda\|_1 \right\}\right)\, .
\end{align*}
Now observe that both measures $\delta_{X_{n+1}}$ and $\eta_n$ are supported on degrees at most $d_{\mathrm{max}}+n$, where $d_{\mathrm{max}}$ is the maximum degree in the initial graph $G_0$. Hence both $\|\delta_{X_{n+1}}Q^\Lambda \|_1$ and $\|\eta_n Q^\Lambda \|_1$ are smaller (up to constant factors) than $(n+1)^\gamma$. Finally, we obtain that for $t\in [\tau_n,\tau_{n+1}[$,
\begin{align*}
\sup_{s\geq 0} \| A_t(s)\|_1 &\leq \int_t^{+\infty}\|F(X(u))-F(\overline{X}(u))\|_1\d u\\
&\lesssim \sum_{i\geq n}\gamma_{i+1}(i+1)^{\gamma-1} \lesssim \sum_{i\geq n}(i+1)^{\gamma-2}\lesssim n^{\gamma-1} \underset{n\to +\infty}{\longrightarrow}0\, ,
\end{align*}
establishing the first part of~\eqref{eq:A_t-B_t}. Moving on to $B_t$, if $t\in [\tau_{n}, \tau_{n+1})$ and $t+s\in [\tau_{m}, \tau_{m+1})$ with $m\geq n$, then
\[
B_t(s)=(\tau_{n+1}-t)\varepsilon'_{n+1}+\sum_{i=n+1}^{m-1}\gamma_{i+1}\varepsilon'_{i+1}+ (t+s-\tau_m)\varepsilon'_{m+1}\, .
\]
Noticing that $\gamma_{i+1}\varepsilon'_{i+1}=\tfrac{1}{i+1}\varepsilon_{i+1}$, we have
\[
\|B_t(s)\|_1\leq \frac{1}{n+1}\|\varepsilon_{n+1}\|_1 +\frac{1}{m+1}\|\varepsilon_{m+1}\|_1 +\left\| \sum_{i=n+1}^{m-1}\frac{1}{i+1}\varepsilon_{i+1}\right\|_1\, .
\]
For the first two terms, we can argue as above: 
\[
\|\varepsilon_{n+1}\|_1\lesssim 1+\max_{\Lambda\in\{A, B\}}\|\eta_n Q^\Lambda\|_1\lesssim (n+1)^\gamma\, ,
\]
so that $\tfrac{1}{n+1}\|\varepsilon_{n+1}\|_1 \to 0$ as $n\to +\infty$. For the last term, we show that
\begin{equation}\label{eq:convergence-noise}
    \sup_{m\geq n}\left\| \sum_{i=n}^m \frac{1}{i+1}\varepsilon_{i+1}\right\|_1\, \underset{n\to +\infty}{\longrightarrow}\, 0\quad \text{a.s.}
\end{equation}
The arguments are quite similar to the ones that were used in the proof of Lemma~\ref{lem:compact-closure} to prove that the quantity $L$ was finite almost surely. The proof can be found in Section~\ref{sec:proof-convergence-noise}. This yields $\sup_{s\geq 0} \| B_t(s)\|_1  \to 0$ as $t\to +\infty$, establishing~\eqref{eq:A_t-B_t}. 

Hence, if $X^\star$ is a limit point of $\{\Theta^t(X)\}_{t\geq 0}$, then it satisfies
\[
X^\star=L_F(X^\star)\, ,
\]
that is, $X^\star$ is solution to~\eqref{eq:dynamicsystem-1}.
\end{proof}

\section{Existence and uniqueness of the fixed point}\label{sec:fixed-point-unique}

For $u \in [0, 1]$, let $Q_u = u Q^A + (1-u)Q^B$. Then, for $\eta \in \mathcal{P}(E)$, we have $F(\eta) = \eta Q_{s(\eta)} - \eta Q_{s(\eta)}(E) \eta$, with $s(\eta)$ defined in~\eqref{eq::def-s}. In this section, we prove that there exists a unique probability measure $\eta^\star$ on $E$ such that $F(\eta^\star) = 0$. To this end, we show that for all $u \in [0, 1]$, there exists a unique probability measure $\pi_u$ on $E$ such that $\pi_u$ is a left eigenvector of $Q_u$, then it follows from a straightforward calculation that the associated eigenvalue $\lambda(u) > 0$. Furthermore, we observe that if a probability measure $\eta$ on $E$ satisfies $F(\eta) = 0$, then $\eta$ is a left eigenvector of $Q_{s(\eta)}$ associated with the eigenvalue $\lambda_\eta := \eta Q_{s(\eta)}(E) > 0$. Since $s(\eta) \in [0, 1]$, there exists some $u$ such that $\eta = \pi_u$, and this parameter satisfies $s(\pi_u) = u$. The uniqueness of a solution to the equation $F(\eta) = 0$ follows from the uniqueness of a solution to the equation $s(\pi_u) = u$. We denote this unique solution by $u^\star$ (or $s^\star$ since $s(\pi_{u^\star}) = u^\star$), so that $\eta^\star = \pi_{u^\star}$. Finally, we present an alternative proof of this result based on Schauder's fixed-point theorem, which yields a different, yet still implicit, characterization of the solution $\eta^\star$.

\subsection{Existence and uniqueness of a eigenmeasure of $Q_u$}

Let us first recall the definition of $Q_u$. For every $u \in [0, 1]$, $Q_u$ is the convex combination of $Q^A$ and $Q^B$ defined as $Q_u := u Q^A + (1-u)Q^B$. More explicitly, $Q_u$ can be written as the sum of the following two infinite matrices:

\[
\begin{pmatrix}
up & (1-u)(1-p) & 0 & 0 & \cdots \\
up & (1-u)(1-p) & 0 & 0 & \cdots \\
up & (1-u)(1-p) & 0 & 0 & \cdots \\
up & (1-u)(1-p) & 0 & 0 & \cdots \\
\vdots & \vdots & \vdots & \vdots & \ddots
\end{pmatrix}
\]
and
\[
\scriptsize{
\renewcommand{\arraystretch}{1.3}
\begin{blockarray}{c*{10}{c}}
 & \lab{(1, A)} & \lab{(1, B)} & \lab{(2, A)} & \lab{(2, B)} & \lab{\cdots}
 & \lab{(k, A)} & \lab{(k, B)} & \lab{(k\!+\!1, A)} & \lab{(k\!+\!1, B)} & \lab{\cdots} \\
 & \downarrow & \downarrow & \downarrow & \downarrow & &
   \downarrow & \downarrow & \downarrow & \downarrow & \\
\begin{block}{c(*{10}{c})}
\lab{(1, A)}\;\rightarrow & -u & 0 & u2^{\gamma} & 0 & \cdots & 0 & 0 & 0 & 0 & \cdots \\
\lab{(1, B)}\;\rightarrow & 0 & -(1-u) & 0 & (1-u)2^{\gamma} & \cdots & 0 & 0 & 0 & 0 & \cdots \\
\lab{(2, A)}\;\rightarrow & 0 & 0 & -u2^{\gamma} & 0 & \cdots & 0 & 0 & 0 & 0 & \cdots \\
\lab{(2, B)}\;\rightarrow & 0 & 0 & 0 & -(1-u)2^{\gamma} & \cdots & 0 & 0 & 0 & 0 & \cdots \\
\lab{\vdots} & \vdots & \vdots & \vdots & \vdots & \ddots & \vdots & \vdots & \vdots & \vdots & \\
\lab{(k-1, A)}\;\rightarrow & 0 & 0 & 0 & 0 & \cdots & uk^{\gamma} & 0 & 0 & 0 & \cdots \\
\lab{(k-1, B)}\;\rightarrow & 0 & 0 & 0 & 0 & \cdots & 0 & (1-u)k^{\gamma} & 0 & 0 & \cdots \\
\lab{(k, A)}\;\rightarrow & 0 & 0 & 0 & 0 & \cdots & -uk^{\gamma} & 0 & u(k\!+\!1)^{\gamma} & 0 & \cdots \\
\lab{(k, B)}\;\rightarrow & 0 & 0 & 0 & 0 & \cdots & 0 & -(1-u)k^{\gamma} & 0 & (1-u)(k+1)^{\gamma}  & \cdots \\
\lab{\vdots} & \vdots & \vdots & \vdots & \vdots & & \vdots & \vdots & \vdots & \vdots & \ddots \\
\end{block}
\end{blockarray}
}\, .
\]

\begin{lemma}\label{lem:eigenmeasure}
For every $u\in[0, 1]$, there exists a unique probability measure on $E$ that is a eigenmeasure of $Q_u$. We denote this measure by $\pi_u$ and we have 
\[
\pi_u(k, A) = p\displaystyle\prod\limits_{j=1}^{k} \cfrac{u j^\gamma}{\lambda(u)+u j^\gamma} \quad \text{and} \quad \pi_u(k, B) = (1-p)\displaystyle\prod\limits_{j=1}^{k} \cfrac{(1-u)j^\gamma}{\lambda(u)+(1-u)j^\gamma}\, ,
\]
where $\lambda(u)$ denotes the unique real number satisfying
\[
p\displaystyle\sum\limits_{k=1}^{\infty} \displaystyle\prod\limits_{j=1}^{k} \cfrac{u j^\gamma}{\lambda(u)+u j^\gamma} + (1-p)\displaystyle\sum\limits_{k=1}^{\infty}\displaystyle\prod\limits_{j=1}^{k} \cfrac{(1-u)j^\gamma}{\lambda(u)+(1-u)j^\gamma} = 1\, .
\]
In other words, $(\lambda(u), \pi_u)$ is the unique eigenpair of $Q_u$ such that $\lambda(u) \in \R$ and $\pi_u \in \mathcal{P}(E)$.
Moreover, it turns out that $\lambda(u)$ is equal to $\pi_u Q_u(E)$ and is strictly positive.
\end{lemma}

\begin{proof}[Proof]
Suppose that there exist $\eta\in\mathcal{P}(E)$ and $\lambda\in\mathbb{R}$ such that $\eta Q_u = \lambda \eta$. The equality $\eta Q_u = \lambda \eta$ is equivalent, componentwise, to the following system:
\[
\begin{cases}
    \eta Q_u(1, A) = up-u\eta(1, A) &= \lambda \eta(1, A) \\
    \eta Q_u(1, B) = (1-u)(1-p)-(1-u)\eta(1, B) &= \lambda \eta(1, B) \\
    \eta Q_u(k, A) = u k^\gamma (\eta(k-1, A)-\eta(k, A)) &= \lambda \eta(k, A), \quad \forall k \geq 2 \\
    \eta Q_u(k, B) = (1-u)k^\gamma (\eta(k-1, B)-\eta(k, B)) &= \lambda \eta(k, B), \quad \forall k \geq 2 \, .
\end{cases}
\]
Thus, $(\lambda, \eta)$ satisfies the following relations:
\[
\begin{cases}
    \eta(1, A) = p \;\cfrac{u}{\lambda +u} \\
    \eta(1, B) = (1-p) \;\cfrac{1-u}{\lambda + 1-u}
\end{cases}
\]
and
\[
\forall k \geq 2, \quad
\begin{cases}
    \eta(k, A) = \cfrac{u k^\gamma}{\lambda +u k^\gamma} \; \eta(k-1, A) \\
    \eta(k, B) = \cfrac{(1-u) k^\gamma}{\lambda +(1-u) k^\gamma} \; \eta(k-1, B) \, .
\end{cases} 
\]
By induction, we have
\begin{equation}\label{eq:stationnarymeasure}
\forall k \geq 1, \quad
\begin{cases}
    \eta(k, A) = p \displaystyle\prod\limits_{j=1}^{k} \cfrac{u j^\gamma}{\lambda +u j^\gamma} \\
    \eta(k, B) = (1-p)\displaystyle\prod\limits_{j=1}^{k} \cfrac{(1-u) j^\gamma}{\lambda +(1-u) j^\gamma} \, .
\end{cases} 
\end{equation}
Note that, for $\eta$ to be a probability measure, all the terms appearing in~\eqref{eq:stationnarymeasure} must be positive and sum to $1$. It follows that $\lambda$ must be strictly positive. Let $\chi$ be the function defined on $(0, \infty)$ by
\[
\chi(x) = p \displaystyle\sum\limits_{k=1}^{\infty} \displaystyle\prod\limits_{j=1}^{k} \cfrac{u j^\gamma}{x +u j^\gamma} + (1-p) \displaystyle\sum\limits_{k=1}^{\infty} \displaystyle\prod\limits_{j=1}^{k} \cfrac{(1-u) j^\gamma}{x +(1-u) j^\gamma} \, . 
\]
Since $\eta$ must be a probability measure, such a pair $(\eta, \lambda)$ exists if and only if the equation $\chi(x) = 1$ admits a solution $x > 0$. The function $\chi$ is clearly strictly decreasing on $(0, \infty)$. Moreover, by standard arguments, we can show that $\chi:(0, \infty)\to(0, \infty)$ is continuous and thus bijective. Hence, the equation $\chi(x) = 1$ admits a unique solution, which we denote by $\lambda(u)$. Defining $\pi_u$ as in~\eqref{eq:stationnarymeasure} with $\lambda = \lambda(u)$, we have thus established the existence and uniqueness of a left eigenpair $(\lambda(u), \pi_u) \in \R \times \mathcal{P}(E)$ of $Q_u$.\\
The fact that $\pi_u Q_u = \lambda(u) \pi_u$, combined with the fact that $\pi_u$ is a probability measure, immediately implies that $\lambda(u) = \pi_u Q_u(E)$.
\end{proof}

\subsection{Existence and uniqueness of a solution to $F(\eta) = 0$}

\begin{prop}\label{prop:stationnarysolution}
Let $F$ be the function defined on $\mathcal{P}(E)$ by $F(\eta) = \eta Q_{s(\eta)} - \eta Q_{s(\eta)}(E) \eta$. There exists a unique probability measure $\eta\in\mathcal{P}(E)$ such that $F(\eta)=0$.
\end{prop}
To prove this proposition, we use the following lemma, whose proof is provided at the end of the paper.
\begin{lemma}\label{lem:fixedpoint}
Let $s$ be the function defined on $\mathcal{P}(E)$ by
\[
s(\eta) = \cfrac{r_A \eta Q^B(E)}{r_B \eta Q^A(E) + r_A \eta Q^B(E)}\, ,
\]
where we recall the notation $r_\Lambda = p_\Lambda (1-q_\Lambda) + p_{\overline{\Lambda}} q_{\overline{\Lambda}}$. For every $u \in [0, 1]$, recall that $(\lambda(u),\pi_u)$ is the unique eigenpair in $\R \times \mathcal{P}(E)$ of $Q_u$ . There exists a unique $u \in [0, 1]$ such that $s(\pi_u) = u$. We denote by $s^\star$ this unique solution and we set $s_A^\star = s^\star$ and $s_B^\star = 1-s^\star$.
\end{lemma}

\begin{proof}[Proof of Proposition~\ref{prop:stationnarysolution}]
Suppose that there exists a probability measure $\eta$ on $E$ such that $F(\eta) = 0$. Thus, $\eta$ is a eigenmeasure of $Q_{s(\eta)}$ associated with the eigenvalue $\lambda_\eta := \eta Q_{s(\eta)}(E)$. Since $u:=s(\eta) \in [0, 1]$, Lemma~\ref{lem:eigenmeasure} ensures that $\eta = \pi_{s(\eta)}=\pi_u$. Moreover, $\pi_{s(\eta)}$ satisfies $\pi_{s(\eta)} Q_{s(\eta)} = \lambda_\eta \pi_{s(\eta)}$. And it follows from this equality that $s(\eta) = s(\pi_{s(\eta)})$, which, recalling the notation $u :=s(\eta)$, can be rewritten as $u = s(\pi_u)$. The existence and uniqueness of a solution to this equation, established by Lemma~\ref{lem:fixedpoint}, completes the proof.
\end{proof}

\begin{remark}
Note that, in the single-community case, the stationary measure $\eta^\star$ is given by, for all $k \in \N^*$,
\[
\eta^\star(k) = \prod_{j=1}^{k} \cfrac{j^\gamma}{\lambda^\star + j^\gamma}\, ,
\]
where $\lambda^\star$ is the unique solution $\lambda$ to the following equation
\[
\sum_{k=1}^{\infty} \prod_{j=1}^{k} \cfrac{j^\gamma}{\lambda + j^\gamma} = 1 \,.
\]
This parameter $\lambda^\star$ corresponds to the Malthusian parameter defined in~\cite{rudas2007random}.
\end{remark}

\subsection{Alternative proof}

Let us give an alternative proof of the existence and uniqueness of $\eta^\star\in\cP(E)$ such that $F(\eta^\star)=0$.

\begin{proof}[Alternative proof of Proposition~\ref{prop:stationnarysolution}]
Consider the operator $R$ defined by: for all $\mu\in\cP(E)$,
\[
R_\mu=s(\mu)(Q^A-\mu Q^A(E)I)+(1-s(\mu))(Q^B-\mu Q^B(E)I)\, ,
\]
and define the application
\begin{align*}
    \Psi: \quad & \cP(E)\longrightarrow \cP(E)\\
    & \mu \longmapsto \Psi(\mu)\, ,
\end{align*}
where $\Psi(\mu)$ is the unique probability measure $\pi$ such that $\pi R_\mu=0$. One can check that $\Psi(\mu)$ is given by: for all $k\geq 1$ and $\Lambda\in\{A, B\}$,
\[
\Psi(\mu)(k, \Lambda)=C^{-1}p_\Lambda\prod_{j=1}^k\frac{j^\gamma}{j^\gamma+\frac{\alpha^\Lambda_\mu}{r_\Lambda}}\, ,
\]
where 
\[
\alpha^\Lambda_\mu=\mu Q^\Lambda(E)=p_\Lambda +\sum_{k\geq 1}\mu(k, \Lambda)\left((k+1)^\gamma-k^\gamma\right)\, ,
\]
and where $C>0$ is the normalizing constant:
\[
C=\sum_{\Lambda\in\{A, B\}}\sum_{k\geq 1} p_\Lambda\prod_{j=1}^k\frac{j^\gamma}{j^\gamma+\frac{\alpha^\Lambda_\mu}{r_\Lambda}}\, \cdot 
\]
We observe that $\Psi$ is continuous and that, since $\alpha^\Lambda_\mu$ is always between $p_\Lambda$ and $p_\Lambda+2^\gamma-1$, all the probability measures $\Psi(\mu)$ admit a first moment that is uniformly bounded by some constant $M>0$. Hence $\Psi$ takes values in a compact convex subset of the Banach space $\ell^1(E)$ and by Schauder's fixed point theorem, it admits a fixed point. In other words, there exists a probability measure $\eta^\star\in\cP(E)$ such that $\eta^\star R(\eta^\star)=F(\eta^\star)=0$. 

Let us now show that this measure is unique. A solution $\eta^\star$ satisfies, for all $k\geq 1$, $\Lambda\in\{A, B\}$,
\[
\eta^\star(k, \Lambda)=p_\Lambda\prod_{j=1}^k\frac{j^\gamma}{j^\gamma+\frac{\alpha^\Lambda_{\eta^\star}}{r_\Lambda}}\, \cdot 
\]
In particular, we have
\[
\alpha^\Lambda_{\eta^\star}=p_\Lambda+\sum_{k\geq 1}p_\Lambda\prod_{j=1}^k\frac{j^\gamma}{j^\gamma+\frac{\alpha^\Lambda_{\eta^\star}}{r_\Lambda}}\left((k+1)^\gamma-k^\gamma\right)\, .
\]
Since the function $x\mapsto \sum_{k\geq 1}\prod_{j=1}^k\frac{j^\gamma}{j^\gamma+\frac{x}{r_\Lambda}}\left((k+1)^\gamma-k^\gamma\right)$ is strictly decreasing on $\R_+^*$, the values $\alpha^A_{\eta^\star}$ and $\alpha^B_{\eta^\star}$ are uniquely determined, hence $\eta^\star$ is unique.
\end{proof}

\begin{remark}\label{rem:link-lambda-alpha}
The two proofs yield two expressions for the stationary measure $\eta^\star$, one parametrized by $\lambda^\star$ and $s^\star = s(\eta^\star)$ and the other by $\alpha_{\eta^\star}^\Lambda$, for $\Lambda \in \{A, B\}$. There is, of course, a natural relationship between these three parameters $\lambda^\star$, $s(\eta^\star)$ and $\alpha_{\eta^\star}^\Lambda$. Since we have
\begin{align*}
    \lambda^\star &= \eta^\star Q_{s(\eta^\star)}(E) = s^\star_A \alpha_{\eta^\star}^A +  s^\star_B \alpha_{\eta^\star}^B \\
    &= \frac{r_A \alpha_{\eta^\star}^B \alpha_{\eta^\star}^A}{r_B \alpha_{\eta^\star}^A + r_A \alpha_{\eta^\star}^B} + \frac{r_B \alpha_{\eta^\star}^A \alpha_{\eta^\star}^B}{r_B \alpha_{\eta^\star}^A + r_A \alpha_{\eta^\star}^B} \\
    &= \frac{\alpha_{\eta^\star}^A \alpha_{\eta^\star}^B}{r_B \alpha_{\eta^\star}^A + r_A \alpha_{\eta^\star}^B} \qquad (\text{since } r_A + r_B = 1),
\end{align*}
it follows that
\[
\frac{\lambda^\star}{s^\star_\Lambda} = \frac{\alpha_{\eta^\star}^\Lambda}{r_\Lambda}\, ,
\]
for $\Lambda \in \{A, B\}$.
\end{remark}

\section{Convergence in the one-community case}\label{sec:one-community}

In the case of one single community, the operator $F$ that drives the equation $\partial_t \eta_t=F(\eta_t)$ is given by
\begin{equation}\label{eq:dynamics-one-community}
F(\eta)=\eta Q - \alpha_\eta\eta\, ,
\end{equation}
where $Q$ is given by
\[
Q=
\begin{pmatrix}
1 & 0 & 0 & \dots & \dots \\
1 & 0 & 0 & \dots & \dots \\
\vdots & \vdots & \vdots & \dots &\dots \\
1 & 0 & 0 & \dots &\dots \\
\vdots & \vdots & \vdots & \dots &\dots \\
\end{pmatrix}
+ \begin{pmatrix}
-1^\gamma & 2^\gamma & 0 & \dots & \dots & \\
0 & -2^\gamma & 3^\gamma & 0 & \dots &\\
0 & 0 & \ddots & \ddots& &\\
0 & \dots & & -k^\gamma & (k+1)^\gamma &\\
0 & \dots & \dots &   & \ddots & \ddots
\end{pmatrix}
\]
and where
\[
\alpha_\eta=\eta Q(E):=1+ \sum_{k\geq 1}\eta_k\left((k+1)^\gamma -k^\gamma\right)\, .
\]
One approach to study Equation~\eqref{eq:dynamics-one-community} is to follow~\cite{mailler2020stochastic} and to view it as the equation describing the law of a conditioned Markov process. Noticing that the sum on each row of the matrix $Q$ is always less that $2^\gamma$, the matrix $Q-2^\gamma I$ is a sub-Markovian infinitesimal generator: upon adding an additional state $\partial$ with the appropriate jump rates, it becomes Markovian. The function $F(\eta)=\eta (Q-2^\gamma I)-\eta(Q-2^\gamma I)(E)\eta$ is then interpreted as the one that governs the distribution of a Markov process on $E\cup \{\partial\}$ conditioned to stay on $E$, and one may use results such as~\cite{ChampagnatVillemonais2016,ChampagnatVillemonais2023} to show convergence to a unique \emph{quasi-stationary distribution}. In our setting, the problem of this approach is that it is not clear how to extend it to two communities. Instead, our approach consists in trying to write $F(\eta)$ as $F(\eta)=\eta L_\eta$ with $L_\eta$ having negative entries on the diagonal, positive entries off the diagonal, and zero sum on each row. Note that $Q-\alpha_\eta I$ is not far from doing the job but the sum on each row does not equal zero. One may rather consider the following small perturbation of $Q-\alpha_\eta I$ and define $L_\eta$ as
\[
L_\eta =
\begin{pmatrix}
{-2^\gamma}
& {2^\gamma}
& & & & \\[3pt]
\alpha_\eta+2^\gamma-3^\gamma
& {-(\alpha_\eta+2^\gamma)}
& {3^\gamma}
& & & \\
\vdots & \qquad \ddots & \qquad \ddots &  & \qquad (0) & \\
\vdots & & \qquad \ddots & \qquad \ddots &  & \\
\vdots & & & \qquad \ddots & \qquad \ddots & \\
\alpha_\eta+k^\gamma-(k+1)^\gamma
& \qquad (0) & &
& {-(\alpha_\eta+k^\gamma)}
& {(k+1)^\gamma} & \cdots \\
\vdots &  &  & & \qquad \ddots & \quad \ddots
\end{pmatrix}\, .
\]
In words, the first column of the matrix $L_\eta $ is given by the sequence $(\alpha_\eta +k^\gamma-(k+1)^\gamma)$ except for the first entry which is $-2^\gamma$, the diagonal is given by the sequence $(-\alpha_\eta -k^\gamma)$ (except for the first entry which is $-2^\gamma$), and the top-diagonal is given by the sequence $((k+1)^\gamma)$. Note that, for all $k\geq 2$, we have
\[
\alpha_\eta+k^\gamma-(k+1)^\gamma\geq \theta:= 1+2^\gamma-3^\gamma>0\, . 
\]
The operator $L_\eta$ corresponds to the infinitesimal generator of a non-linear Markov jump process\footnote{See~\cite{chaintron2022propagationI}[Section A.4.3] for a nice introduction to non-linear Markov processes.} on the positive integers which jumps from $k$ to $k+1$ at rate $(k+1)^\gamma$, and from $k$ to $1$ at rate $\alpha_\eta +k^\gamma -(k+1)^\gamma$. We observe that only the jump rates associated with returns to state $1$ depend on the law, whereas the jump rates from $k$ to $k+1$ depend only on $k$. For $\gamma\in\{0, 1\}$, the non-linearity completely disappears from $L_\eta$: whatever the current state $k\geq 2$, the process jumps to $1$ at rate $1$. Convergence can be easily established in those two extreme cases, \emph{e.g.}\ by standard coupling arguments for (linear) Markov processes.

\subsection{Construction of the process $(X_t)_{t \geq 0}$} 

Let $(\mu_t)_{t \geq 0}$ be a continuous trajectory taking values in $\mathcal{P}(E)$. The equation 
\begin{equation}\label{eq:eq_law_time_inhomogeneous_Markov_process}
    \partial \eta_t = \eta_t L_{\mu_t}
\end{equation}
is then the evolution equation for the law of a time-inhomogeneous Markov process whose infinitesimal generator at time $t$ is given by $L_{\mu_t}$. In this section, we construct a process $(X_t)_{t \geq 0}$ whose law $(\eta_t)_{t \geq 0}$ is a solution to Equation~\eqref{eq:eq_law_time_inhomogeneous_Markov_process}. To this end, we show that, for every bounded measurable function $f_0$ on $E$, the function $f_t := P_{0, t}^{\mu} f_0$ associated with the semi-group $(P_{0, t}^{\mu})_{t \geq 0}$ is a solution to $\partial_t f_t = L_{\mu_t} f_t$. Note that we use the same notation $L$ for the action of an operator on functions and on measures, the dual action being distinguished only by whether the operator acts on the left (for functions) or on the right (for measures). Then in Section~\ref{sec:thecoupling-1}, we construct a coupling that allows us to establish, in Section~\ref{sec:existence-uniqueness1}, the existence and uniqueness of a solution to the equation $\partial_t \eta_t = \eta_t L_{\eta_t}$.

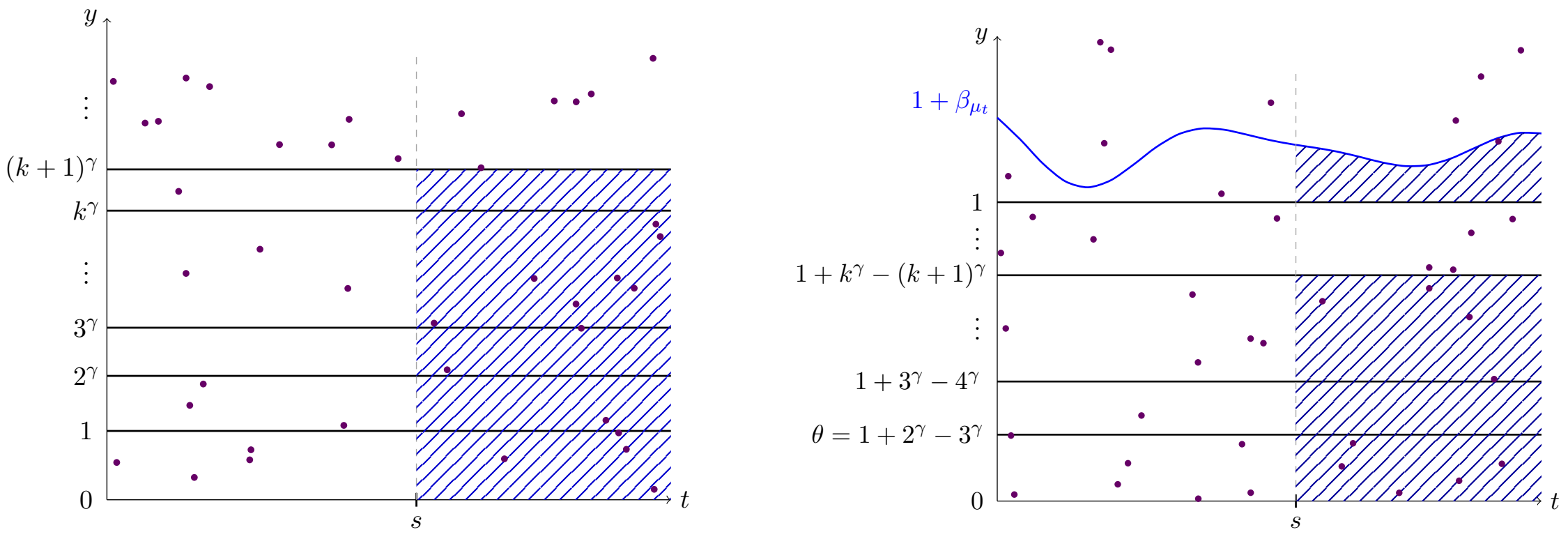
\begin{figure}[ht] 
\centering

\makebox[\textwidth][c]{%
\begin{minipage}{1\textwidth} 
\centering

\begin{subfigure}{0.45\textwidth} 
\centering 
\resizebox{\textwidth}{!}{%
\begin{tikzpicture}[domain=0:8.2]
\draw[->] (0,0) -- (8.2,0) node[right] {$t$};
\draw[->] (0,0) -- (0,7) node[left] {$y$};
\node at (-0.3,0) {$0$} ;
\draw[thick]  (0,1) -- (8.2,1) ;
\node at (-0.3,1) {$1$} ;
\draw[thick]  (0,1.8) -- (8.2,1.8) ;
\node at (-0.3,1.8) {$2^\gamma$} ;
\draw[thick]  (0,2.5) -- (8.2,2.5) ; 
\node at (-0.3,2.5) {$3^\gamma$} ;
\node at (-0.3,3.4) {$\vdots$} ;
\draw[thick]  (0,4.2) -- (8.2,4.2) ;
\node at (-0.3,4.2) {$k^\gamma$} ;
\draw[thick]  (0,4.8) -- (8.2,4.8) ;
\node at (-0.8,4.8) {$(k+1)^\gamma$} ;
\node at (-0.3,5.8) {$\vdots$} ;

    \draw[thick] (4.5, 0.1) -- (4.5, -0.1) node[below] {$s$};
    \draw[dashed, gray!70] (4.5, 0) -- (4.5, {6.5});

    \fill[pattern={Lines[angle=45, distance=3pt]}, pattern color=blue!80!black] 
        plot[smooth, domain=4.5:8.2] (\x, {4.8}) -- 
        plot[smooth, domain=8.2:4.5] (\x, {0}) -- cycle;
        
    \pgfmathsetseed{65} 
    \foreach \i in {1,...,40} {
    \pgfmathsetmacro{\x}{8.2 * rnd} 
    \pgfmathsetmacro{\y}{6.5 * rnd} 
    \fill[violet!80!black] (\x,\y) circle (1.4pt);} 
\end{tikzpicture} }
\caption{Point cloud $\Pi_{\text{up}}$} \label{fig:point-cloud-up} 
\end{subfigure}
\hfill 
\begin{subfigure}{0.5\textwidth} 
\centering 
\resizebox{\textwidth}{!}{%
\begin{tikzpicture}[domain=0:8.2]

\draw[->] (0,0) -- (8.2,0) node[right] {$t$};
\draw[->] (0,0) -- (0,7) node[left] {$y$};

\node at (-0.3,0) {$0$} ;
\draw[thick]  (0,1) -- (8.2,1) ;
\node at (-1.5,1) {$\theta=1+2^\gamma-3^\gamma$} ;
\draw[thick]  (0,1.8) -- (8.2,1.8) ;
\node at (-1.2,1.8) {$1+3^\gamma-4^\gamma$} ;
\node at (-0.3,2.7) {$\vdots$} ;
\draw[thick]  (0,3.4) -- (8.2,3.4) ;   
\node at (-1.6,3.4) {$1+k^\gamma-(k+1)^\gamma$} ;
\node at (-0.3,4.05) {$\vdots$} ;
\draw[thick]  (0,4.5) -- (8.2,4.5) ;
\node at (-0.3,4.5) {$1$} ;
\node[blue] at (-0.7,6) {$1+\beta_{\mu_t}$};

\draw[thick, blue] plot[smooth] (\x, {5.3 + 0.3*cos(\x*80+40) + 0.2*cos(\x*100 + 50) + 0.15*cos(\x*130+40)});

    \draw[thick] (4.5, 0.1) -- (4.5, -0.1) node[below] {$s$};
    \draw[dashed, gray!70] (4.5, 0) -- (4.5, {6.5});

\fill[pattern={Lines[angle=45, distance=3pt]}, pattern color=blue!60!black] 
        (4.5,0) rectangle (8.2,3.4);

\fill[pattern={Lines[angle=45, distance=3pt]}, pattern color=blue!60!black] 
    plot[smooth, domain=4.5:8.2] (\x, {5.3 + 0.3*cos(\x*80+40) + 0.2*cos(\x*100 + 50) + 0.15*cos(\x*130+40)}) -- 
    (8.2, 4.5) -- (4.5, 4.5) -- cycle;

\pgfmathsetseed{46} 
\foreach \i in {1,...,40} {
    \pgfmathsetmacro{\x}{8.2 * rnd} 
    \pgfmathsetmacro{\y}{7 * rnd} 
    \fill[violet!80!black] (\x,\y) circle (1.4pt);
}

\end{tikzpicture} } 
\caption{Point cloud $\Pi_{\text{down}}$} \label{fig:point-cloud-down} 
\end{subfigure} 
\caption{Representation of the two point clouds $\Pi_{\text{up}}$ and $\Pi_{\text{down}}$. The hatched regions indicate where to look for the point corresponding to the next jump when $X$ is in state $k$ at time $s$.}
\label{fig:point-clouds}
\end{minipage}%
}
\end{figure}

Let us consider two independent Poisson point processes $\Pi_{\text{up}}$ and $\Pi_{\text{down}}$ on $\R_+ \times \R_+$ with intensity $\mathrm{d}t \times \mathrm{d}y$. The first process governs the upward jumps, while the second governs the returns to $1$. Let us first focus on the upward jump rates. As illustrated in Figure~\ref{fig:point-cloud-up}, which corresponds to the point cloud $\Pi_{\text{up}}$, for each $k \in \mathbb{N}$, we draw the half-line defined by the equation $y = k^\gamma$ on $\R_+$. We thus partition the cone $\mathbb{R}_+ \times \mathbb{R}_+$ into horizontal strips of respective heights $1-0, \:2^\gamma-1, \: 3^\gamma-2^\gamma, ... , \: (k+1)^\gamma - k^\gamma$, \textit{etc}. For the point cloud $\Pi_{\text{down}}$, as illustrated in Figure~\ref{fig:point-cloud-down}, we partition the cone in the same way, by drawing the horizontal lines $y = 1 + k^\gamma - (k+1)^\gamma$, for each $k\in\N^*$, and then the curve $t \mapsto \beta_{\eta_t}$, for $ t \in\R_+$, where
\begin{equation}\label{eq:def-beta}
\beta_\eta:= \sum_{k\geq 1}\eta_k\left((k+1)^\gamma -k^\gamma\right)=\alpha_\eta-1\, .
\end{equation}

Now, let $\eta_0 \in \mathcal{P}(E)$ be the initial condition for Equation~\eqref{eq:eq_law_time_inhomogeneous_Markov_process}. We then define the random variable $X_0$ to be distributed according to $\eta_0$. We define the stopping time $\tau_1$ as the random variable corresponding to the first jump time of the process. More precisely, if $X_0$ is at state $k$ for some $k \in \N^*$, $\tau_1$ is given by the point of minimal abscissa $t>0$ among the points of $\Pi = \Pi_{\text{up}} \cup \Pi_{\text{down}}$ that lie, for the point cloud $\Pi_{\text{up}}$, between the lines $y=0$ and $y = (k+1)^\gamma$, and, for the point cloud $\Pi_{\text{down}}$, either between the lines $y=0$ and $y=1+k^\gamma-(k+1)^\gamma$, or between the line $y=1$ and the curve $y=1+\beta_{\eta_t}$. Note that $\tau_1$ is well defined, since these regions have finite area and, almost surely, contain only finitely many points, with no two points having the same abscissa. If the point belongs to $\Pi_{\text{up}}$, then $X_{\tau_1}=X_0 +1 = k+1$, whereas if it belongs to $\Pi_{\text{down}}$, then $X_{\tau_1}= 1$. Moreover, for every $t\in[0,\tau_1)$ --- note that $\tau_1$ is strictly positive almost surely ---, we have $X_t=X_0=k$. We proceed in the same way for the second jump, and so on. Figure~\ref{fig:point-clouds} shows the regions in which to look for the point with minimal abscissa when the process $X$ is in state $k$ at time $s$, for the two point clouds. Thus, by induction, we define a sequence of stopping times $(\tau_n)_{n \in \N^*}$ that is strictly increasing almost surely and hence the process $X$ for all $t \in [0,\sup_{n \geq 0} \tau_n)$ --- where we set $\tau_0 = 0$. In fact, this construction allows us to define the process $X$ on the whole of $\R_+$, since $X$ is non-explosive, as established in the following lemma.

\begin{lemma}\label{lem:non-explosive}
The process $X$ constructed above is non-explosive, that is,
\[
\tau_\infty := \sup_{n \geq 0} \tau_n = \infty \quad \mathrm{a.s.}
\]
\end{lemma}

\begin{proof}
To begin with, observe that the process $(X_t)_{t \geq 0}$ is stochastically dominated by the homogeneous Markov process $(Y_t)_{t \geq 0}$, with $Y_0 := X_0$, that jumps from $k$ to $k+1$ at rate $(k+1)^\gamma$, for each $k \in \N^*$, and from $k \geq 2$ to $1$ at rate $c_\gamma$, where the constant $c_\gamma$ is defined by
\[
c_\gamma = \begin{cases}
    1 + 2^\gamma - 3^\gamma & \quad \text{if } \; 0 \leq \gamma <1 \\
    1 & \quad \text{if} \; \gamma = 1 \, .
\end{cases}
\]
Note that this stochastic domination implies that we can construct a coupling $(X_t, Y_t)_{t \geq 0}$ with $X_0 = Y_0$ (we keep the same notation $X$ and $Y$ for the coupled processes) such that whenever $X$ jumps upward, $Y$ does so as well, and whenever $Y$ returns to state $1$, so does $X$. In particular, this coupling ensures that $Y$ is always greater than or equal to $X$, which is an equivalent characterization of the fact that the law of $Y$ stochastically dominates that of $X$. For each $n\geq 1$, let $\tau_n^X$ and $\tau_n^Y$ denote the $n$-th jump times of the processes $X$ and $Y$, respectively. The previous coupling ensures that $\tau_n^X \geq \tau_n^Y$ for every $n\geq 1$, almost surely. We therefore have almost surely
\[
\sup_{n \geq 1} \tau_n^X \geq \sup_{n \geq 1} \tau_n^Y \, .
\]
Thus, it suffices to show that $\sup_{n \geq 1} \tau_n^Y = +\infty$ almost surely to prove that the process $X$ is non-explosive. However, by standard results on continuous-time Markov jump processes on a countable state space, showing that the process $Y$ is recurrent implies that $\lim_{n\to\infty}\tau_n^Y=+\infty$ almost surely. Moreover, for the model associated with the process $Y$, it is actually sufficient to show that state $1$ is recurrent, since the probability of reaching state $k \geq 2$ starting from state $1$ is strictly positive for every $k \geq 2$. Let us now show that state $1$ is recurrent for the process $Y$. Conditional on $Y_n = k$ ($k \geq 2$), the probability that $Y_{n+1}=k+1$ is equal to 
\[
\frac{(k+1)^\gamma}{(k+1)^\gamma + c_\gamma} \, .
\] 
To show that state $1$ is recurrent, it suffices to consider the process started from state $1$ and prove that its first return time to state $1$, denoted by $T_1^+$, is almost surely finite. We have
\[
\PP_1(T_1^+ = \infty) = \lim_{n \to \infty} \PP_1(T_1^+ > n) \, ,
\]
where we denote by $\PP_1$ the probability measure under which $Y$ is a continuous-time Markov jump process with the transition rates specified above and initial condition $Y_0 = 1$ almost surely. Hence, we have
\begin{align*}
\PP_1(T_1^+ > n) &= \PP_1(Y_1 = 2, Y_2 = 3, ... , Y_n = n+1) \\
&= \prod_{k=1}^{n} \frac{(k+1)^\gamma}{(k+1)^\gamma +c_\gamma} \\
&\leq \exp\left(-\sum_{k=1}^{n} \cfrac{c_\gamma}{(k+1)^\gamma + c_\gamma}\right) \underset{n\to +\infty}{\longrightarrow} 0 \, ,
\end{align*}
as $\gamma \in [0, 1]$ and $c_\gamma > 0$. Thus, $T_1^+$ is finite almost surely, which concludes the proof.
\end{proof}

Finally, it is straightforward to verify that the law of the process $X$ constructed above is a solution to Equation~\eqref{eq:eq_law_time_inhomogeneous_Markov_process}.

\subsection{The coupling}\label{sec:thecoupling-1} 
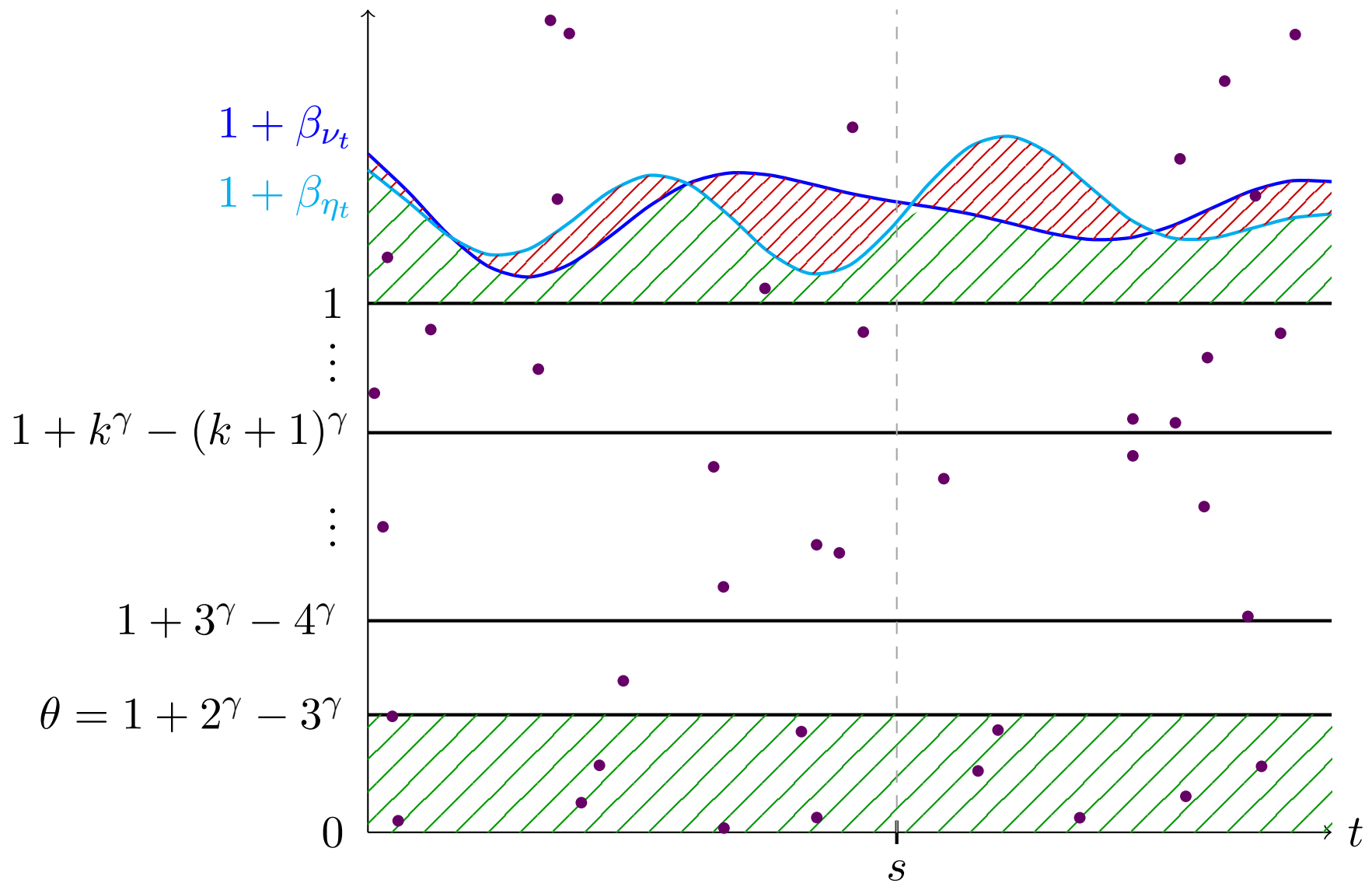
\begin{figure}
    \centering

\begin{tikzpicture}[domain=0:8.2]

\draw[->] (0,0) -- (8.2,0) node[right] {$t$};
\draw[->] (0,0) -- (0,7) ;

\node at (-0.3,0) {$0$} ;
\draw[thick]  (0,1) -- (8.2,1) ;
\node at (-1.5,1) {$\theta=1+2^\gamma-3^\gamma$} ;
\draw[thick]  (0,1.8) -- (8.2,1.8) ;
\node at (-1.2,1.8) {$1+3^\gamma-4^\gamma$} ;
\node at (-0.3,2.7) {$\vdots$} ;
\draw[thick]  (0,3.4) -- (8.2,3.4) ;    
\node at (-1.6,3.4) {$1+k^\gamma-(k+1)^\gamma$} ;
\node at (-0.3,4.1) {$\vdots$} ;
\draw[thick]  (0,4.5) -- (8.2,4.5) ;
\node at (-0.3,4.5) {$1$} ;
\node[blue] at (-0.7,6) {$1+\beta_{\nu_t}$};
\node[cyan] at (-0.7,5.4) {$1+\beta_{\eta_t}$};
    
\fill[pattern={Lines[angle=45, distance=5.5pt]}, pattern color=green!60!black] 
        (0,0) rectangle (8.2,1);

\fill[pattern={Lines[angle=45, distance=5.5pt]}, pattern color=green!60!black] 
    plot[smooth, domain=0:8.2] (\x, {min(5.3 + 0.3*cos(\x*80+40) + 0.2*cos(\x*100 + 50) + 0.15*cos(\x*130+40), 5.3 + 0.2*cos(\x*60+20) + 0.3*cos(\x*110+100) + 0.2*cos(\x*140)}) -- 
    (8.2, 4.5) -- (0, 4.5) -- cycle;

\draw[thick, blue] plot[smooth] (\x, {5.3 + 0.3*cos(\x*80+40) + 0.2*cos(\x*100 + 50) + 0.15*cos(\x*130+40)});
\draw[thick, cyan] plot[smooth] (\x, {5.3 + 0.2*cos(\x*60+20) + 0.3*cos(\x*110+100) + 0.2*cos(\x*140)});

    \fill[pattern={Lines[angle=45, distance=3pt]}, pattern color=red!80!black] 
        plot[smooth, domain=0:8.2] (\x, {5.3 + 0.3*cos(\x*80+40) + 0.2*cos(\x*100 + 50) + 0.15*cos(\x*130+40)}) -- 
        plot[smooth, domain=8.2:0] (\x, {5.3 + 0.2*cos(\x*60+20) + 0.3*cos(\x*110+100) + 0.2*cos(\x*140)}) -- cycle;

    \draw[thick] (4.5, 0.1) -- (4.5, -0.1) node[below] {$s$};
    \draw[dashed, gray!70] (4.5, 0) -- (4.5, {7});

\pgfmathsetseed{46} 
\foreach \i in {1,...,40} {
    \pgfmathsetmacro{\x}{8.2 * rnd} 
    \pgfmathsetmacro{\y}{7 * rnd} 
    \fill[violet!80!black] (\x,\y) circle (1.4pt);
}

\end{tikzpicture}
    \caption{Point cloud $\Pi_{\text{down}}$ associated to downward jumps to $1$. We refer to the finely hatched region as the ``red region'' and to the more sparsely hatched region as the ``green region''.}
    \label{fig:coupling-down}
\end{figure}

Now, consider two trajectories $(\mu^1_t)_{t \geq 0}$ and $(\mu^2_t)_{t \geq 0}$ taking values in $\mathcal{P}(E)$. In this section, we construct a coupling of two stochastic processes $(X_t)_{t \geq 0}$ and $(Y_t)_{t \geq 0}$ whose laws at time $t$ are respectively given by $\eta^1_t$, where $(\eta^1_t)_{t \geq 0}$ solves 
\begin{equation}\label{eq:eq-mu-1}
    \partial \eta_t = \eta_t L_{\mu^1_t} \, ,
\end{equation}
and $\eta^2_t$, where $(\eta^2_t)_{t \geq 0}$ solves 
\begin{equation}\label{eq:eq-mu-2}
    \partial \eta_t = \eta_t L_{\mu^2_t}\, .
\end{equation}
Let $\eta_0^1$ and $\eta_0^2$ be two probability measures on $E$, corresponding to the initial conditions for Equations~\eqref{eq:eq-mu-1} and~\eqref{eq:eq-mu-2}, respectively. Let $(X_0, Y_0)$ be an optimal coupling of $\eta_0^1$ and $\eta_0^2$. As before, let $\Pi_{\text{up}}$ and $\Pi_{\text{down}}$ be two independent Poisson point processes on $\R_+ \times \R_+$ with intensity $\mathrm{d}t \times \mathrm{d}y$. For both processes $X$ and $Y$, we use the same Poisson point processes: the first governs the upward jumps, while the second governs the returns to $1$. The only difference is that, for the jumps back to state $1$, the curve $y = 1 + \beta_{\mu_t^1}$ is used to construct $X$, whereas the curve $y = 1 + \beta_{\mu_t^2}$ is used to construct $Y$, as illustrated in Figure~\ref{fig:coupling-down}, with $\beta_\eta$ defined in~\eqref{eq:def-beta}. As a consequence, when the two processes are in the same state $k \in \N^*$, they can separate only when the first point appears in the ``red region'' (see Figure~\ref{fig:coupling-down} for the definitions of the ``red'' and ``green regions''). We will exploit this key property of the coupling both to establish the existence and uniqueness of a solution to the dynamical system $\partial_t \eta_t = \eta_t L_{\eta_t}$ and to prove that this solution converges to the stationary state as $t$ tends to infinity.

\subsection{Existence and uniqueness of the solution}\label{sec:existence-uniqueness1} 

We can now establish the existence and uniqueness result stated earlier by combining the coupling constructed above with a contraction argument.

\begin{prop}\label{prop:existence-uniqueness-one-community}
For all measure $\eta_0\in\mathcal{P}(E)$, there exists a unique solution to $\partial_t \eta_t= \eta_t L_{\eta_t}$.
\end{prop}

\begin{proof}[Proof]
Let $\mu_1=(\mu^1_t)_{t\geq 0}$ and $\mu_2=(\mu^2_t)_{t\geq 0}$ be two trajectories of measures in $\mathcal{P}(E)$, and let $\nu_0$ be a measure in $\mathcal{P}(E)$. Consider the process $(X_t, Y_t)_{t\geq 0}$ as defined in Section~\ref{sec:thecoupling-1}. Since $X_0=Y_0$, we have
\begin{align*}
    \| \nu_0 P_{0, t}^{\mu_1}-\nu_0 P_{0, t}^{\mu_2}\|_{\mathrm{TV}} &\leq \P(\tau<t)\, ,
\end{align*}
where $\tau$ corresponds to the arrival time of the first point in the ``red zone'' (see Figure~\ref{fig:coupling-down} for the definition of the ``red region''). Indeed, when the two processes are in the same state at time $t = 0$, the coupling is constructed so that they first separate when one process jumps back to state $1$ while the other remains unchanged. This corresponds precisely to a jump triggered by a point in the ``red region''. By setting $\beta_s^i : = \beta_{\mu_s^i}$, we have
\begin{align*}
\P(\tau<t) &\leq 1-\exp\left(-\int_0^t|\beta_s^1-\beta_s^2|\mathrm{d} s\right)\\
&\leq 1-\exp\left(-(2^\gamma-1)\int_0^t\|\mu_s^1-\mu_s^2\|_{\mathrm{TV}}\mathrm{d} s \right)\, ,
\end{align*}
where we used the characterization of the total-variation distance 
$$
\|\mu-\nu\|_{\mathrm{TV}}=\sup_{f,\; 0\leq f\leq 1}|\E_\mu f-\E_\nu f|\, ,
$$
and the fact that $\beta_\mu =\E_\mu f$, with, for all $k\geq 1$, $f(k)=(k+1)^\gamma -k^\gamma\leq 2^\gamma-1$. Using the inequality $1-\e^{-x}\leq x$, we obtain
\begin{align*}
    \| \nu_0 P_{0, t}^{\mu_1}-\nu_0 P_{0, t}^{\mu_2}\|_{\mathrm{TV}} &\leq (2^\gamma-1)\int_0^t\|\mu_s^1-\mu_s^2\|_{\mathrm{TV}}\mathrm{d} s\, .
\end{align*}
Now, setting $t_0=\frac{1}{2(2^\gamma-1)}$, we see that
\[
  \sup_{s\in [0, t_0]}\| \nu_0 P_{0, s}^{\mu_1}-\nu_0 P_{0, s}^{\mu_2}\|_{\mathrm{TV}} \leq \frac{1}{2}\sup_{s\in [0, t_0]}\|\mu_s^1-\mu_s^2\|_{\mathrm{TV}}\, .
\]
In particular, the application $\Psi: (\mu_t)_{t\in [0, t_0]} \mapsto (\nu_0 P_{0, t}^\mu)_{t\in [0, t_0]}$ is a contraction of the Banach space $L^\infty([0, t_0], \mathcal{P}(E))$. Since this space is complete, the Banach fixed-point theorem ensures that $\Psi$ admits a unique fixed point, which is a solution of $\partial_t \eta_t= \eta_t L_{\eta_t}$ on $[0, t_0]$. Then there exists a unique solution on
$[nt_0, (n + 1)t_0]$ for all $n\geq 0$, and thus on $\R_+$.
\end{proof}

\begin{remark}
Rewriting $F(\eta)$ as $\eta L_\eta$ allows us to interpret the equation $\partial_t \eta_t = F(\eta_t)$ as a nonlinear Fokker–Planck-type evolution equation (on a countable state space), that is, as the deterministic equation governing the evolution of the probability distribution. The construction of the process $(X_t)_{t \geq 0}$ we have presented naturally leads to a McKean–Vlasov-type equation, which describes the evolution of a stochastic process whose dynamics depend on its own probability law, as governed by the associated nonlinear Fokker–Planck equation. The two-part survey by Chaintron and Diez provides a very accessible and engaging introduction to topics related to propagation of chaos and, consequently, to McKean–Vlasov equations (see~\cite{chaintron2022propagationI} and~\cite{chaintron2022propagationII}).
\end{remark}

\subsection{Convergence of the solution}\label{sec:proof-convergence-1}

\begin{prop}\label{prop:coupling-one-community}
Let $\gamma \in [0, 1]$ satisfy $\beta^\star > 3^\gamma - 2$ or $\gamma=1$. Let $(\eta_t)_{t\geq 0}$ be the unique solution to $\partial_t \eta_t= \eta_t L_{\eta_t}$ on $\R_+$, and $\eta^\star$ be unique solution to $\eta L_\eta=0$. We have
\[
\|\eta_t-\eta^\star\|_{\mathrm{TV}}\, \underset{t\to \infty}{\longrightarrow}\, 0\, .
\]
\end{prop}
\begin{proof}[Proof]
Let $(X_t, Y_t)_{t\geq 0}$ be a coupling of $(\eta_t)_{t\geq 0}$ and $\eta^\star$, generated as in Section~\ref{sec:thecoupling-1}, with $(X_0, Y_0)$ distributed according to an optimal coupling of $\eta_0$ and $\eta^\star$. Let 
\[
\tau=\inf\{s\geq 0, X_s=Y_s\}\, .
\]
This random variable has an atom at $0$ since $\P(\tau=0)=1-\|\eta_0-\eta^\star\|_{\mathrm{TV}}$, and is absolutely continuous with respect to the Lebesgue measure $\lambda$ on $\R_+^*$. Letting $f_\tau$ be the density of $\tau$ with respect to $\mu=\delta_0+\lambda$, we have
\begin{align*}
    \P(X_t\neq Y_t)&=\P(\forall s\in [0, t],\, X_s\neq Y_s)+\int_{0}^t\P(X_t\neq Y_t\,|\, \tau=s)f_\tau(s)\mathrm{d} \mu(s)\, .
\end{align*}
On the one hand, by setting $\theta := 1+2^\gamma-3^\gamma$, we have
\[
\P(\forall s\in [0, t],\, X_s\neq Y_s)\leq \e^{-\theta t}\|\eta_0-\eta^\star\|_{\mathrm{TV}}\, ,
\]
because the presence of a point at time $s$ between the lines $y=0$ and $y=\theta$ brings both processes back to state $1$, so that they coincide at time $s$. On the other hand, conditionally on $\{\tau=s\}$, if $X_t\neq Y_t$, then there is a point in the ``red zone'' between $s$ and $t$ which is not followed by a point in the ``green zone'' (see Figure~\ref{fig:coupling-down} for the definitions of the ``green'' and ``red regions''). Let $\tau'$ be the random variable corresponding to the time of the last point in the ``red zone'' between $s$ and $t$ if there is one, and to $-\infty$ otherwise. To compute the density of $\tau'$ observe that, for $u\in [s, t]$, denoting by $\mathcal{P}(\lambda)$ a Poisson random variable with parameter $\lambda$,
\begin{align*}
    \P(\tau'\leq u)&= \P\left(\mathcal{P}\left(\int_u^t |\beta_x -\beta^\star|\mathrm{d} x\right)=0\right) = \exp\left(-\int_u^t |\beta_x -\beta^\star|\mathrm{d} x\right)\, ,
\end{align*}
where $\beta_x=\beta_{\eta_x}$. Hence $\tau'$ has a density on $[s, t]$ given by, for all $u\in [s, t]$,
\[
f_{\tau'}(u)=|\beta_u -\beta^\star|\exp\left(-\int_u^t |\beta_x -\beta^\star|\mathrm{d} x\right)\, .
\]
Noticing that the probability that there is no point in the ``green zone'' between $u$ and $t$ is equal to $\e^{-\int_u^t (\theta +\min\{\beta_x, \beta^\star\})\mathrm{d} x}$, this yields
\begin{align*}
    \P(X_t\neq Y_t\,|\, \tau=s)&\leq \int_s^t f_{\tau'}(u)\e^{-\int_u^t (\theta +\min\{\beta_x, \beta^\star\})\mathrm{d} x} \mathrm{d} u\\
    &\leq \int_s^t |\beta_u -\beta^\star|\e^{-\int_u^t (\theta +\max\{\beta_x,\beta^\star\})\mathrm{d}x}\mathrm{d} u\\
    &\leq \int_s^t |\beta_u -\beta^\star|\e^{-\int_u^t (\theta +\beta^\star)\mathrm{d}x}\mathrm{d} u\\
    &\leq \int_s^t |\beta_u -\beta^\star|\e^{-(\theta +\beta^\star)(t-u)}\mathrm{d} u \, .
\end{align*}
Notice that the right-hand side above is maximized when $s = 0$, so that we get
\[
\P(X_t\neq Y_t)\leq \e^{-\theta t}\|\eta_0-\eta^\star\|_{\mathrm{TV}}+\int_0^t |\beta_u -\beta^\star|\e^{-(\theta +\beta^\star)(t-u)}\mathrm{d} u \, .
\]
Using that $|\beta_u -\beta^\star|\leq (2^\gamma-1)\|\eta_u-\eta^\star\|_{\mathrm{TV}}$, and setting $g(t)=\e^{(\theta +\beta^\star) t}\|\eta_t-\eta^\star\|_{\mathrm{TV}}$, we have
\[
    g(t)\leq \e^{\beta^\star t} g(0) +(2^\gamma-1)\int_0^t g(u)\mathrm{d} u \, .
\]
By Gronwall's lemma, we have
\begin{align*}
g(t)&\leq g(0)\left( \e^{\beta^\star t}  +(2^\gamma-1)\int_0^t \e^{\beta^\star u} \e^{(2^\gamma-1)(t-u)} \mathrm{d} u\right)\\
&= g(0)\left( \e^{\beta^\star t}  +(2^\gamma-1)\frac{\e^{(2^\gamma-1)t}-\e^{\beta^\star t}}{2^\gamma-1-\beta^\star}\right) \, .
\end{align*}
Therefore,
\begin{equation}\label{eq:inequality-contraction-total-variation}
    \|\eta_t -\eta^\star\|_{\mathrm{TV}} \leq \|\eta_0 -\eta^\star\|_{\mathrm{TV}}\left( \e^{-\theta t}+\frac{2^\gamma-1}{2^\gamma-1-\beta^\star}\left(\e^{-(\beta^\star +2-3^\gamma)t} - \e^{-\theta t}\right)\right) \, .
\end{equation}
When $\gamma<1$, then $\theta > 0$. In this case, if $\beta^\star+2-3^\gamma > 0$, then the claimed convergence follows. When $\gamma = 1$, the ``red region'' in Figure~\ref{fig:coupling-down} is empty, as $\alpha_\eta = 1 + \beta_\eta = 2$ for every $\eta \in \mathcal{P}(E)$ and thus we have
\begin{equation*}
    \|\eta_t -\eta^\star\|_{\mathrm{TV}} \leq \|\eta_0 -\eta^\star\|_{\mathrm{TV}} \, \e^{- t} \, . \qedhere
\end{equation*}
\end{proof}

\begin{lemma}\label{lem:equiv-inequality}
For all $\gamma \in [0, 1]$, the following equivalence holds
\[
\beta^\star >3^\gamma-2\, \iff \sum_{k\geq 1} \prod_{j=1}^k\frac{j^\gamma}{j^\gamma +3^\gamma-1}> 1\, .
\]
\end{lemma}
\begin{proof}[Proof]
First, notice that $\eta^\star L_{\eta^\star}=0$ is equivalent to
\[
\begin{cases}
& \eta^\star(1)=\frac{1}{1+\alpha^\star}\\
& \forall k\geq 2,\; \eta^\star(k)=\frac{k^\gamma}{k^\gamma +\alpha^\star} \eta^\star(k-1)\, ,
\end{cases}
\]
where $\alpha^\star=1+\beta^\star$. Hence, for all $k\geq 1$, we have
\[
\eta^\star(k)=\prod_{j=1}^k \frac{j^\gamma}{j^\gamma +1+\beta^\star}\, ,
\]
and $\beta^\star$ can be defined as the solution of
\[
\beta=\sum_{k\geq 1} \prod_{j=1}^k \frac{j^\gamma}{j^\gamma +1+\beta}\left((k+1)^\gamma-k^\gamma\right)\, .
\]
Since the function $\beta \mapsto \sum_{k\geq 1} \prod_{j=1}^k \frac{j^\gamma}{j^\gamma +1+\beta}\left((k+1)^\gamma-k^\gamma\right)$ is decreasing, showing that $\beta^\star > 3^\gamma -2$ amounts to showing that 
\[
\sum_{k\geq 1} \prod_{j=1}^k \frac{j^\gamma}{j^\gamma +3^\gamma-1}\left((k+1)^\gamma-k^\gamma\right)>3^\gamma -2\, .
\]
After simplifications, this is equivalent to 
\begin{equation}\label{eq:inequality-coupling}
    \sum_{k\geq 1} \prod_{j=1}^k\frac{j^\gamma}{j^\gamma +3^\gamma-1}> 1\, ,
\end{equation}
which concludes the proof.
\end{proof}

\begin{figure}[!ht]
    \centering
    \scalebox{0.8}{
    \begin{tikzpicture}
        \begin{axis}[
            width=12cm,
            height=8cm,
            xlabel={$\gamma$},
            ylabel={$\varphi(\gamma)$},
            xmin=0,
            xmax=1,
            ymin=0,
            ymax=12,
            axis lines=left,
            grid=major,
            tick align=outside,
        ]
            \addplot[
                black,
                thick,
                no marks
            ]
            table[
                x=x,
                y=f,
                col sep=space
            ] {fonction.dat};
            \addplot[
                red,
                thick,
                domain=0:1,
                samples=2
            ] {1};
            \node[
                red,
                anchor=south west
            ] at (axis cs:0.01,1.1) {$y=1$};
        \end{axis}
    \end{tikzpicture}
    }
    \caption{Graph of
    $\displaystyle
    \varphi(\gamma)=\sum_{k=1}^{+\infty}
    \prod_{j=1}^{k}
    \frac{j^\gamma}{j^\gamma+3^\gamma-1}
    $
    on $(0, 1]$.}
    \label{fig:function-inequality}
\end{figure}
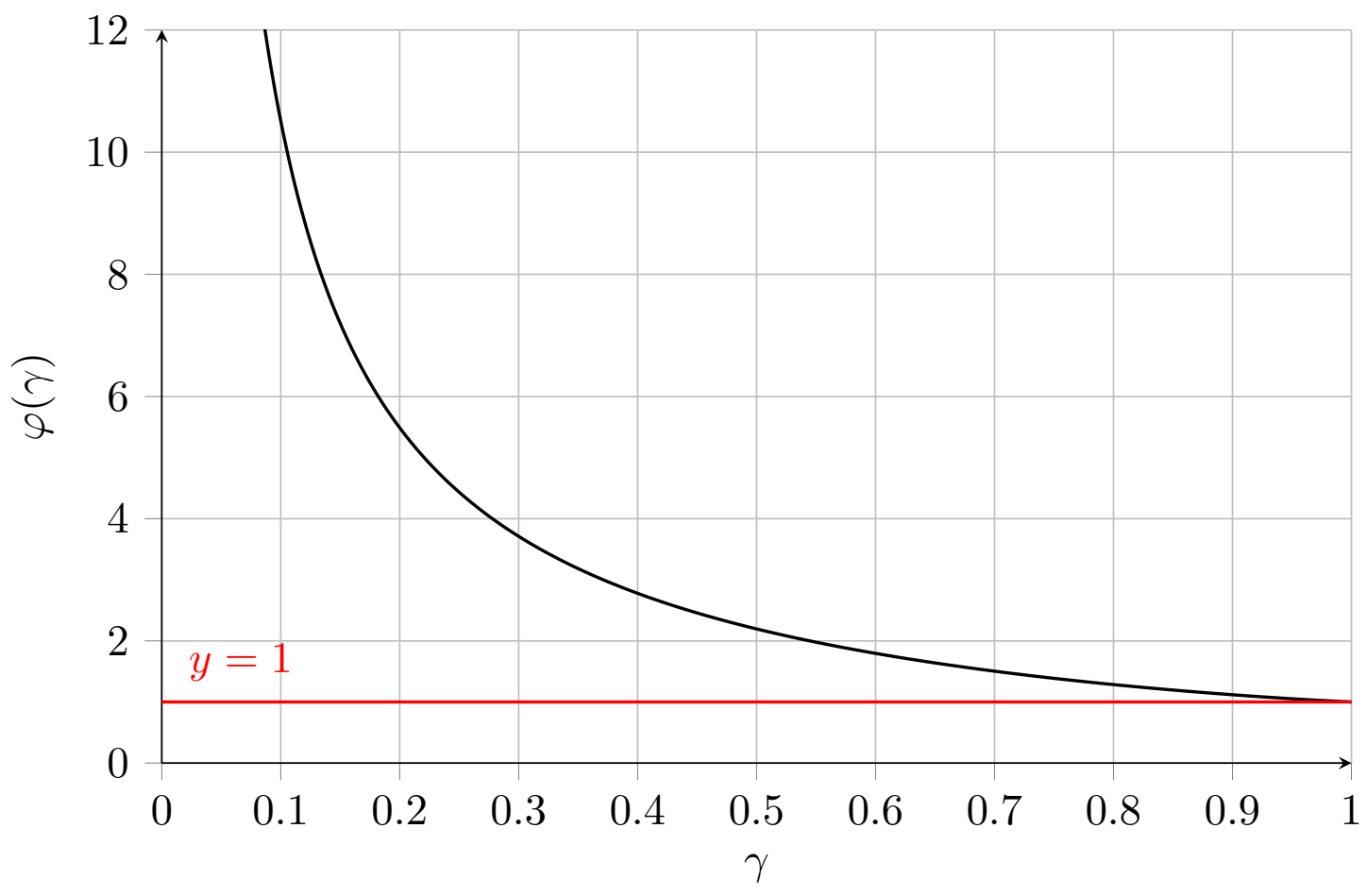

The preceding lemma allows us to obtain an equivalent reformulation of Proposition~\ref{prop:coupling-one-community}. Since proving Inequality~\eqref{eq:inequality-coupling} is equivalent to showing that $\beta^\star > 3^\gamma -2$, and since $\beta_\eta > 0$ for all $\eta \in \mathcal{P}(E)$, the inequality $\beta^\star > 3^\gamma -2$ holds for every $\gamma \in [0, 1]$ such that $3^\gamma-2 \leq 0$, that is, for every $\gamma \in [0, \frac{\ln(2)}{\ln(3)}]$. We have not been able to find an analytical proof of this inequality for all $\gamma \in [0, 1]$. However, the graph of the function $\varphi$ (see Figure~\ref{fig:function-inequality}), together with the fact that its curve lies locally above the line $y=1$ in a left-hand neighborhood of $\gamma =1$ (see Remark~\ref{rem:derivative-at-1}), strongly suggests that the inequality holds throughout $[0, 1]$. Nevertheless, we strongly believe that a numerical proof of this inequality can be obtained by using the mean value inequality to construct a sufficiently fine mesh $(x_i)_{1\leq i\leq N}$ of $[0, 1]$, consisting of $N$ points, such that the inequality holds at each mesh point $x_i$. This can be established by numerically approximating the infinite sum by a finite sum and obtaining a suitable bound on the remainder. The conclusion then follows by controlling the variation of $\varphi$ between consecutive mesh points using the mean value inequality.
Thus, this approach also requires establishing an upper bound on $\lvert \varphi'\rvert$, either globally on $[0, 1]$ or locally on each subinterval $[x_i, x_{i+1}]$. Since such a proof would provide little insight and be of limited interest, we have chosen not to include pages of calculations and numerical approximations.
We therefore leave the search for an elegant proof as an open question.

\begin{remark}\label{rem:derivative-at-1}
Let $\varphi$ be the function defined on $[0, 1]$ by
\[
\varphi(\gamma) = \sum_{k=1}^{\infty} \prod_{j=1}^{k} \frac{j^\gamma}{j^\gamma + 3^\gamma - 1} \, .
\]
Note that $\varphi(0) = \infty$ and $\varphi(1) = 1$. Furthermore, a direct computation shows that the derivative of the function $\varphi$ is equal to
\[
\varphi'(\gamma) = \sum_{k=1}^{\infty} \left(\prod_{j=1}^{k} \cfrac{j^\gamma}{j^\gamma + 3^\gamma - 1}\right) .  \sum_{i=1}^{k} \frac{3^\gamma (\ln i - \ln 3) - \ln i}{i^\gamma + 3^\gamma - 1}\, ,
\]
and that
\begin{align*}
    \varphi'(1) &= \sum_{k=1}^{\infty} \left(\prod_{j=1}^{k} \cfrac{j}{j+2} \right) \sum_{i=1}^{k} \frac{2 \ln i - 3 \ln 3}{i+2} \\
    &= \sum_{i=1}^{\infty} \frac{2 \ln i - 3 \ln 3}{i+2} \sum_{k=i}^{\infty} 2 \left(\frac{1}{k+1} - \frac{1}{k+2}\right) \\
    &= 4 \, \sum_{i=1}^{\infty} \frac{\ln i}{(i+1)(i+2)} - 3 \ln 3 \, .
\end{align*}
A comparison between the series and the corresponding integral shows that
\[
\sum_{i=1}^{\infty} \frac{\ln i}{(i+1)(i+2)} \leq \int_{1}^{\infty} \frac{\ln x}{x(x+1)} \mathrm{d}x \, .
\]
Finally, standard integration techniques show that the latter integral is equal to $\frac{\pi^2}{12}$. Hence, we obtain
\[
\varphi'(1) \leq \frac{\pi^2}{3} - 3 \ln 3 <  -0,00597 \, .
\]
The above numerical approximation allows us to conclude that the derivative of the function $\varphi$ at $\gamma = 1$ is strictly negative. Consequently, there exists $\varepsilon>0$ such that Inequality~\eqref{eq:inequality-coupling} holds for all $\gamma \in [1-\varepsilon, 1]$.
\end{remark}

\section{Convergence in the two-community case}\label{sec:convergence-two-communities}

In the same way as in the one-community, we would like to write the dynamics as the one of a non-linear Markov process. To do so, we use the fact that, for all $\eta\in\cP(E)$,
\[
\eta \begin{pmatrix}
p_A\\
p_A\\
\vdots\\
p_A\\
p_A\\
\vdots
\end{pmatrix}
=\eta\begin{pmatrix}
\alpha_\eta^A-2^\gamma +1\\
\alpha_\eta^A\\
\vdots\\
\alpha_\eta^A-(k+1)^\gamma +k^\gamma\\
\alpha_\eta^A\\
\vdots
\end{pmatrix}
\quad \text{ and } \quad 
\eta \begin{pmatrix}
p_B\\
p_B\\
\vdots\\
p_B\\
p_B\\
\vdots
\end{pmatrix}
=\eta\begin{pmatrix}
\alpha_\eta^B\\
\alpha_\eta^B-2^\gamma+1\\
\vdots\\
\alpha_\eta^B\\
\alpha_\eta^B-(k+1)^\gamma +k^\gamma\\
\vdots
\end{pmatrix} \, ,
\]
where for $\Lambda\in\{A, B\}$, we set
\[
\alpha^\Lambda_\eta=\eta Q^\Lambda(E)=p_\Lambda+\sum_{k\geq 1}\eta(k, \Lambda)\left((k+1)^\gamma-k^\gamma\right)\, .
\]
This leads us to define $L^A_\eta$ as
\[\scalebox{0.89}{$
L^A_\eta=
\begin{pmatrix}
\alpha_\eta^A-2^\gamma +1 & 0 &0 &  \dots \\
\alpha_\eta^A & 0 &0 &  \dots \\
\vdots & \vdots &   & \\
\alpha_\eta^A-(k+1)^\gamma +k^\gamma & 0 &0 &  \dots \\
\alpha_\eta^A & 0 &0 &  \dots \\
\vdots & \vdots & & 
\end{pmatrix}
+
\begin{pmatrix}
-1& 0 & 1 & 0 & \dots & & & &\\
0 & -1 & 0 & 1 &0 &\dots & & &\\
\vdots & & \ddots &   & \ddots & & & &\\
  &  &  & -1 & 0& 1 &0 &\dots &\\
 & &  &  & -1 & 0& 1 & 0 &\dots\\
 & & & &  & \ddots& & \ddots &
\end{pmatrix}D^A-\alpha^A_\eta I \, ,
$}
\]
and $L^B_\eta$ as
\[\scalebox{0.89}{$
L^B_\eta=
\begin{pmatrix}
0 & \alpha_\eta^B & 0 & \dots \\
0 & \alpha_\eta^B-2^\gamma +1 & 0 & \dots \\
\vdots & \vdots &  &  \\
0 & \alpha_\eta^B & 0 & \dots \\
0 & \alpha_\eta^B-(k+1)^\gamma +k^\gamma & 0 & \dots \\
\vdots & \vdots & & 
\end{pmatrix}
+
\begin{pmatrix}
-1& 0 & 1 & 0 & \dots & & & &\\
0 & -1 & 0 & 1 &0 &\dots & & &\\
\vdots & & \ddots &   & \ddots & & & &\\
 &  & & -1 & 0& 1 &0 &\dots &\\
 & &  & & -1 & 0& 1 & 0 &\dots\\
 & & & &  & \ddots& & \ddots &
\end{pmatrix}D^B-\alpha^B_\eta I \, .
$}
\]
For $\Lambda\in\{A, B\}$, we have $\eta(Q^\Lambda-\alpha_\eta^\Lambda I)= \eta L_\eta^\Lambda$, and the matrices $L_\eta^\Lambda$ have negative entries on the diagonal and zero sum on each row. However, positivity of the coefficients off the diagonal is not guaranteed: we have to ensure that for all $k\geq 2$, $\alpha_\eta^\Lambda-(k+1)^\gamma +k^\gamma\geq 0$, which amounts to $\alpha_\eta^\Lambda-3^\gamma +2^\gamma\geq 0$. This is why we now make the following assumption.

\begin{assumption}\label{assume:gamma}
The parameter $\gamma\in [0, 1]$ is such that
\[
3^\gamma-2^\gamma \leq p_{\textrm{min}}\, ,
\]
where $p_{\textrm{min}}=\min\{p_A, p_B\}$.
\end{assumption}

\begin{remark}
    It might be the case that this assumption is not necessary in the sense that the dynamics could eventually lead to distributions satisfying Assumption~\ref{assume:gamma}. We were unable to show such a result. While the idea of rewriting the function $F(\eta)$ in the form $\eta L_\eta$ is appealing, since $L_\eta$ can be interpreted as an infinitesimal generator depending on $\eta$, we have nevertheless been unable to establish existence for all values of the parameters $p$ and $\gamma$ in $[0, 1]$.
\end{remark}

We have now rewrote our dynamics as
\begin{equation}\label{eq:new-dynamics-2comm}
\partial_t \eta_t=s(\eta_t)\eta_t L_{\eta_t}^A+(1-s(\eta_t))\eta_t L_{\eta_t}^B\, ,
\end{equation}
which, under Assumption~\ref{assume:gamma}, can be seen as a non-linear Markov process. 

\subsection{Existence and uniqueness of the solution to the non-linear dynamical system}

The next lemma ensures that~\eqref{eq:new-dynamics-2comm} has a unique solution for any initial condition.

\begin{prop}\label{prop:semiflow}
Under Assumption~\ref{assume:gamma}, for any measure $\eta_0\in\cP(E)$, there exists a unique solution to~\eqref{eq:new-dynamics-2comm}.
\end{prop}
\begin{proof}[Proof]
Note that for a fixed sequence $\mu=(\mu_t)_{t\geq 0}\in C^0(\R_+, \cP(E))$, the dynamics 
\[
\partial_t \eta_t=\eta_t\left\{s(\mu_t)(Q^A-\alpha_{\mu_t}^A I)+(1-s(\mu_t))(Q^B-\alpha_{\mu_t}^B I)\right\}
\]
do not preserve probability measures (it preserves positivity of the entries, but not the total mass). However, the dynamics
\begin{equation}\label{eq:dynamics-driven-by-mu}
\partial_t \eta_t=\eta_t\left\{s(\mu_t)L_{\mu_t}^A +(1-s(\mu_t))L_{\mu_t}^B\right\}
\end{equation}
do preserve probability measures. We denote by $(R^\mu_{s, t})_{t\geq s\geq 0}$ the corresponding inhomogeneous Markov semi-group. 

Now let us fix $\nu_0\in\cP(E)$ and define the application
\begin{eqnarray*}
    \Psi:  C^0(\R_+, \cP(E)) &\longrightarrow & C^0(\R_+, \cP(E))\\
    \mu=(\mu_t)_{t\geq 0} & \longmapsto & (\nu_0 R^\mu_{0, t})_{t\geq 0}\, .
\end{eqnarray*}
Let us show that $\Psi$ takes values in the compact convex set $K=K(C)$ defined by
\[
K=\left\{(\nu_t)_{t\geq 0}\in C^0(\R_+, \cP(E)),\, \forall t\geq 0,\, \sum_{x\in E}\nu_t(x)|x|\leq C\right\}\, ,
\]
where if $x=(k, \Lambda)\in E$, then $|x|=k$ (the \emph{height} of $x$), and where $C>0$ will be specified later. In other words, we want to show that, whatever the driving sequence $(\mu_t)$, the expected height along the trajectory is uniformly bounded by the same constant. Then Schauder's fixed point theorem will ensure the existence of a fixed point for $\Psi$, which is a solution to~\eqref{eq:new-dynamics-2comm}. To this end, let $(X^\mu_t)_{t\geq 0}$ be a process whose law is a solution to~\eqref{eq:dynamics-driven-by-mu}. We then define the height process as the first coordinate of the process $(X^\mu_t)_{t \geq 0}$. In the notation introduced above, the height process is therefore given by $(|X^\mu_t|)_{t \geq 0}$. We observe that the height process is always stochastically dominated by the homogeneous Markov process $(H_t)_{t\geq 0}$ on $\N^*$ which starts from the height measure induced by $\nu_0$, and which, when at $k=1$, jumps to $2$ at rate $2^\gamma$, and when at $k\geq 2$, jumps to $k+1$ at rate $(k+1)^\gamma$ and to $1$ at rate 
\begin{equation}\label{eq:def-theta1}
    \theta_1=\inf_{\eta\in\cP(E)}\min_{\Lambda\in\{A, B\}}\left\{s^\Lambda(\eta)\left(\alpha_\eta^\Lambda-3^\gamma+2^\gamma\right)+s^{\overline{\Lambda}}(\eta)\alpha_\eta^{\overline{\Lambda}}\right\}\, .
\end{equation}
For this dominating process, we have
\[
\partial_t \E[H_t]=\E[(H_t+1)^\gamma] -\theta_1(\E[H_t]-1)\leq \left(\E[H_t]+1\right)^\gamma -\theta_1 \E[H_t]+\theta_1\, ,
\]
where we used the Jensen inequality. Letting $h^\star$ be the solution of $(h+1)^\gamma -\theta_1 h+\theta_1$, we have, for all $t\geq 0$,
\[
\E[H_t]\leq \max\{h^\star, h_0\}\, ,
\]
where $h_0=\sum_{x=(k, \Lambda)\in E}\nu_0(x) k$, the initial expected height.
\end{proof}

\begin{remark}
Note that the same proof could have been used for the one-community case. We have chosen instead to present two different proofs: one with a more probabilistic flavor for the one-community case, and the other with a more analytical flavor for the two-community case. The former may have the advantage of being less of a black box. However, this perception may simply reflect our own probabilistic inclinations.

Also note that the two proofs suppose to have written $F(\eta)$ as $\eta L_\eta$, which, in the two-community case, is done under Assumption~\ref{assume:gamma}. One could think of other proofs using \emph{e.g.}\ the Cauchy--Lipschitz Theorem. However, we could not find a norm for which the application $F$ would be Lipschitz and were therefore unable to apply such theorems.
\end{remark}

\subsection{Construction of a process $(X_t)_{t \geq 0}$ whose law solves the dynamical system} Having established the existence and uniqueness of a solution $(\eta_t)_{t \geq 0}$ to Equation~\eqref{eq:new-dynamics-2comm}, we next construct a process $(X_t)_{t \geq 0}$ such that the law of $X_t$ is equal to $\eta_t$, for all $t \in \R_+$, following the approach used in the one-community case. It is worth noting that the time-dependent infinitesimal generator $L_{\eta_t}$ can be interpreted as a convex combination of the two infinitesimal generators $L^A_{\eta_t}$ and $L^B_{\eta_t}$, with weights that depend nonlinearly on $\eta_t$. Equivalently, it can be viewed as the superposition of two dynamics evolving at different rates, namely $s(\eta_t)$ and $1-s(\eta_t)$, respectively. Alternatively, by analogy with a discrete-time evolution, consider the dynamics over a short time interval $[t, t+\varepsilon]$. One can interpret this evolution as follows: at time $t$, a generator is randomly selected by a coin toss and then governs the dynamics throughout the interval $[t, t+\varepsilon]$. The generator $L^A_{\eta_t}$ is selected with probability $s(\eta_t)$, while $L^B_{\eta_t}$ is selected with probability $1-s(\eta_t)$.
Let $\Pi^A_{\text{up}}$, $\Pi^A_{\text{down}}$, $\Pi^B_{\text{up}}$ and $\Pi^B_{\text{down}}$ be four Poisson point processes on $\R_+ \times \R_+$ with intensity $\mathrm{d}t \times \mathrm{d}y$. For $\Lambda \in \{A, B \}$, $\Pi^\Lambda_{\text{up}}$ governs the upward jumps associated with the infinitesimal generator $L^\Lambda$, while $\Pi^\Lambda_{\text{down}}$ governs the returns to state $(1, \Lambda)$. 

Now, let $\eta_0 \in \mathcal{P}(E)$ be the initial condition for Equation~\eqref{eq:new-dynamics-2comm}. We then define the random variable $X_0$ to be distributed according to $\eta_0$. We define the stopping time $\tau_1$ below as the random variable corresponding to the first jump time of the process. A crucial observation is that the form of the infinitesimal generator $L_\eta$ implies that any process whose law solves Equation~\eqref{eq:new-dynamics-2comm} can only jump from state $(k, \Lambda)$ to one of the states $(k+1, \Lambda)$, $(1, A)$, or $(1, B)$. If $X_0=(k, \Lambda)$, with $k\geq 2$, then $\tau_1$ is the point with the smallest abscissa $t > 0$ among the points of $\Pi^\Lambda_{\text{up}} \cup \Pi^\Lambda_{\text{down}} \cup \Pi^{\overline{\Lambda}}_{\text{down}}$ that lie:
\begin{itemize}
    \item \textbf{in the case of $\Pi^\Lambda_{\text{up}}$:} between the curves $y=0$ and $y = s^\Lambda(\eta_t)(k+1)^\gamma$;
    \item \textbf{in the case of $\Pi^\Lambda_{\text{down}}$:} between the curves $y=0$ and $y = s^\Lambda(\eta_t)(\alpha_{\eta_t}^\Lambda + k^\gamma-(k+1)^\gamma)$;
    \item \textbf{in the case of $\Pi^{\overline{\Lambda}}_{\text{down}}$:} between the curves $y=0$ and $y = s^{\overline{\Lambda}}(\eta_t)\alpha_{\eta_t}^{\overline{\Lambda}}$.
\end{itemize}
If $X_0=(1, \Lambda)$, then one only looks at the point with the smallest abscissa $t > 0$ among the points of $\Pi^\Lambda_{\text{up}}$ that lie between the curves $y=0$ and $y=s^\Lambda(\eta_t)2^\gamma$ and those of $\Pi^{\overline{\Lambda}}_{\text{down}}$ that lie between the curves $y= 0$ and $y= s^{\overline{\Lambda}}(\eta_t) \alpha^{\overline{\Lambda}}_{\eta_t}$.

\begin{figure}[ht]
\centering

\begin{minipage}{1\textwidth}
\centering


\begin{subfigure}{0.49\textwidth}
    \centering
    \resizebox{\textwidth}{!}{

\begin{tikzpicture}[domain=0:8.2]

\draw[->] (0,0) -- (8.2,0) node[right] {$t$};
\draw[->] (0,0) -- (0,7) node[left] {$y$};

\node at (-0.3,0) {$0$} ;

\node[violet] at (-0.4,4.7) {$s^{\Lambda}(\eta_t)(\alpha_{\eta_t}^{\Lambda}+k^\gamma - (k+1)^\gamma)$};

\draw[thick, violet] plot[smooth] (\x, {3.5+(0.5+0.2*cos(\x*70+20) + 0.15*cos(\x*120+90) + 0.1*cos(\x*140))*(0.5+0.15*cos(\x*80) + 0.15*cos(\x*100 + 50) + 0.05*cos(\x*120+70))});


\fill[pattern={Lines[angle=45, distance=3pt]}, pattern color=blue!60!black] 
    plot[smooth, domain=3.5:8.2] (\x, {3.5+(0.5+0.2*cos(\x*70+20) + 0.15*cos(\x*120+90) + 0.1*cos(\x*140))*(0.5+0.15*cos(\x*80) + 0.15*cos(\x*100 + 50) + 0.05*cos(\x*120+70))}) -- 
    (8.2, 0) -- (3.5, 0) -- cycle;

    \draw[thick] (3.5, 0.1) -- (3.5, -0.1) node[below] {$s$};
    \draw[dashed, gray!70] (3.5, 0) -- (3.5, {6.5});

\pgfmathsetseed{19} 
\foreach \i in {1,...,40} {
    \pgfmathsetmacro{\x}{8.2 * rnd} 
    \pgfmathsetmacro{\y}{7 * rnd} 
    \fill[violet!80!black] (\x,\y) circle (1.4pt);
}

\end{tikzpicture} 

    }
    \caption{Point cloud $\Pi^{\Lambda}_{\text{down}}$}
    \label{fig:point-cloud-down-2-Lambda}
\end{subfigure}
\hfill
\begin{subfigure}{0.47\textwidth}
    \centering
    \resizebox{\textwidth}{!}{

\begin{tikzpicture}[domain=0:8.2]

\draw[->] (0,0) -- (8.2,0) node[right] {$t$};
\draw[->] (0,0) -- (0,7) node[left] {$y$};

\node at (-0.3,0) {$0$} ;

\node[violet] at (-1,4.8) {$s^{\overline{\Lambda}}(\eta_t)\alpha_{\eta_t}^{\overline{\Lambda}}$};

\draw[thick, violet] plot[smooth] (\x, {4+(0.6+0.2*cos(\x*60+20) + 0.15*cos(\x*110+100) + 0.15*cos(\x*140))*(0.5+0.15*cos(\x*80+10) + 0.1*cos(\x*100 + 30) + 0.05*cos(\x*120+50))});


\fill[pattern={Lines[angle=45, distance=3pt]}, pattern color=blue!60!black] 
    plot[smooth, domain=3.5:8.2] (\x, {4+(0.6+0.2*cos(\x*60+20) + 0.15*cos(\x*110+100) + 0.15*cos(\x*140))*(0.5+0.15*cos(\x*80+10) + 0.1*cos(\x*100 + 30) + 0.05*cos(\x*120+50))}) -- 
    (8.2, 0) -- (3.5, 0) -- cycle;

    \draw[thick] (3.5, 0.1) -- (3.5, -0.1) node[below] {$s$};
    \draw[dashed, gray!70] (3.5, 0) -- (3.5, {6.5});

\pgfmathsetseed{46} 
\foreach \i in {1,...,40} {
    \pgfmathsetmacro{\x}{8.2 * rnd} 
    \pgfmathsetmacro{\y}{7 * rnd} 
    \fill[violet!80!black] (\x,\y) circle (1.4pt);
}

\end{tikzpicture} 
    }
    \caption{Point cloud $\Pi^{\overline{\Lambda}}_{\text{down}}$}
    \label{fig:point-cloud-down-2-overline-Lambda}
\end{subfigure}

\vspace{0.5cm}


\begin{subfigure}{0.5\textwidth}
    \centering
    \resizebox{\textwidth}{!}{

\begin{tikzpicture}[domain=0:8.2]
\draw[->] (0,0) -- (8.2,0) node[right] {$t$};
\draw[->] (0,0) -- (0,7) node[left] {$y$};
\node at (-0.3,0) {$0$} ;

\draw[thick]  (0,5.8) -- (8.2,5.8) ;
\node at (-0.9,5.8) {$(k+1)^\gamma$} ;

\draw[thick, violet] plot[smooth] (\x, {4.4+0.7*(0.5+0.15*cos(\x*100+10) + 0.1*cos(\x*100 + 30) + 0.05*cos(\x*120+50))});
\node[violet] at (-1.4,5) {$s^\Lambda(\eta_t)(k+1)^\gamma$} ;

    \draw[thick] (3.5, 0.1) -- (3.5, -0.1) node[below] {$s$};
    \draw[dashed, gray!70] (3.5, 0) -- (3.5, {6.5});


\fill[pattern={Lines[angle=45, distance=3pt]}, pattern color=blue!60!black] 
    plot[smooth, domain=3.5:8.2] (\x, {4.4+0.7*(0.5+0.15*cos(\x*100+10) + 0.1*cos(\x*100 + 30) + 0.05*cos(\x*120+50))}) -- 
    (8.2, 0) -- (3.5, 0) -- cycle;

    \pgfmathsetseed{65} 
    \foreach \i in {1,...,40} {
    \pgfmathsetmacro{\x}{8.2 * rnd} 
    \pgfmathsetmacro{\y}{6.5 * rnd} 
    \fill[violet!80!black] (\x,\y) circle (1.4pt);} 
\end{tikzpicture}

    }
    \caption{Point cloud $\Pi^\Lambda_{\text{up}}$}
    \label{fig:point-cloud-up-2}
\end{subfigure}


\caption{Representation of the three point clouds $\Pi^\Lambda_{\text{down}}$, $\Pi^{\overline{\Lambda}}_{\text{down}}$ and $\Pi^{\Lambda}_{\text{up}}$. The hatched regions indicate where to look for the point corresponding to the next jump when $X=(k, \Lambda)$ (when $k \geq 2$) at time $s$.}
\label{fig:point-clouds-2}

\end{minipage}
\end{figure}
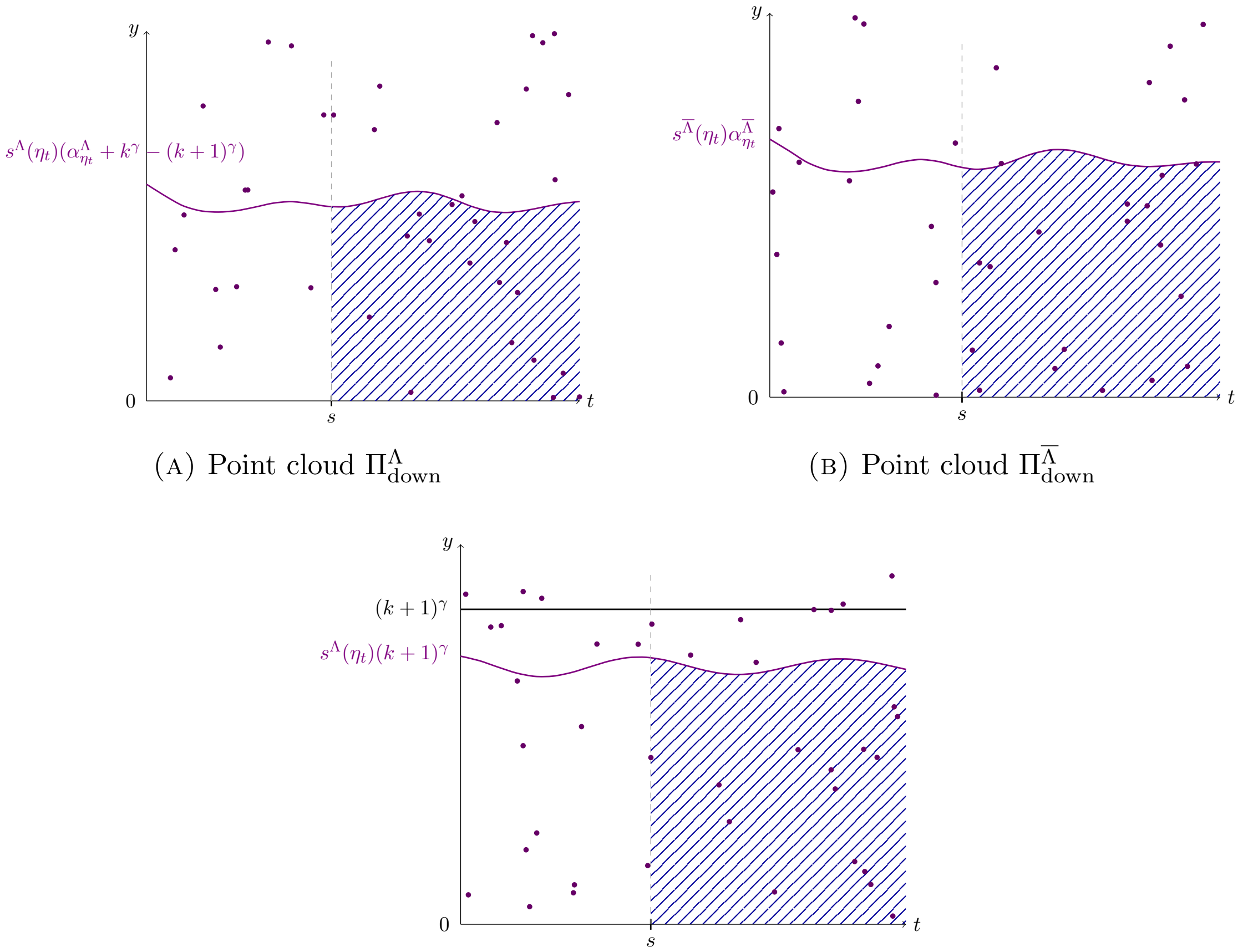

Note that, as in the one-community case, $\tau_1$ is well defined, since these regions have finite area and, almost surely, contain only finitely many points, with no two points having the same abscissa. If the point belongs to $\Pi^\Lambda_{\text{up}}$, then $X_{\tau_1}=(k+1, \Lambda)$, whereas if it belongs to $\Pi^A_{\text{down}}$, then $X_{\tau_1}= (1, A)$ and if it belongs to $\Pi^B_{\text{down}}$, then $X_{\tau_1}= (1, B)$. Moreover, for every $t\in[0, \tau_1)$, we have $X_t=X_0=(k, \Lambda)$. We proceed in the same way for the second jump, and so on. Figure~\ref{fig:point-clouds-2} gives a representation of the regions in which to look for the point with minimal abscissa when the process $X$ is in state $(k, \Lambda)$ (with $k \geq 2$) at time $s$ ($\Pi^{\Lambda}_{\text{down}}$, $\Pi^{\overline{\Lambda}}_{\text{down}}$, and $\Pi^\Lambda_{\text{up}}$ are depicted in Figures~\ref{fig:point-cloud-down-2-Lambda},~\ref{fig:point-cloud-down-2-overline-Lambda} , and~\ref{fig:point-cloud-up-2}, respectively). By induction, we define a sequence of stopping times $(\tau_n)_{n \in \N^*}$ that is, almost surely, strictly increasing and goes to infinity. The fact that $\sup_{n \geq 0} \tau_n = \infty$ almost surely follows from the fact that $\gamma \in [0, 1]$. The jump times and the sequence of states visited thus define a process $(X_t)_{t \geq 0}$ whose law indeed solves Equation~\eqref{eq:new-dynamics-2comm}.

\subsection{The coupling}\label{sec:thecoupling-2} In this section, we construct a coupling of two stochastic processes $(X_t)_{t \geq 0}$ and $(Y_t)_{t \geq 0}$ whose respective laws at time $t$ are $\eta^1_t$ and $\eta^2_t$, two solutions to Equation~\eqref{eq:new-dynamics-2comm}. Let $\eta_0^1$ and $\eta_0^2$ be two probability measures on $E$, corresponding to the initial conditions for this evolution equation. Let $(X_0, Y_0)$ be an optimal coupling of $\eta_0^1$ and $\eta_0^2$. As before, let $\Pi^A_{\text{up}}$, $\Pi^B_{\text{up}}$, $\Pi^A_{\text{down}}$ and $\Pi^B_{\text{down}}$ be four independent Poisson point processes on $\R_+ \times \R_+$ with intensity $\mathrm{d}t \times \mathrm{d}y$. For both processes $X$ and $Y$, we use the same four Poisson point processes: the first two govern the upward jumps along the respective branches $\N^* \times {A}$ and $\N^* \times {B}$ of the graph $E$, while the last two govern the returns to $(1, A)$ and $(1, B)$, respectively. The only difference is that, for each of the four point clouds, the curves drawn on the cone $\R_+ \times \R_+$ depend on the trajectory $(\eta^1_t)_{t \geq 0}$ in the case of $X$ and on $(\eta^2_t)_{t \geq 0}$ in the case of $Y$.

\begin{remark}
Unlike in the one-community case, if the two processes are in the same state $(k, \Lambda)$, with $\Lambda \in \{A, B\}$, at a given time $t$, they can separate either because one jumps to one of the states $(1, A)$ or $(1, B)$ while the other remains unchanged, or because one jumps to state $(k+1, \Lambda)$ while the other does not. In the proof of Proposition~\ref{prop:coupling-one-community}, we made essential use of the fact that the two processes in the constructed coupling could separate only when one of them jumped back to state $1$. In the two-community case, however, the upward jump rates also depend nonlinearly on the law $\eta_t$ at each time $t$, which prevents us from directly adapting the arguments developed in Section~\ref{sec:proof-convergence-1}.
\end{remark}

\subsection{Convergence of the solution}

\begin{prop}\label{prop:coupling-two-community}
Let $(\eta_t)_{t\geq 0}$ be the unique solution to~\eqref{eq:new-dynamics-2comm} on $\R_+$, and $\eta^\star$ be the unique stationary distribution. There exists $0 < \gamma_0 \leq 1$ such that for $\gamma\in [0, \gamma_0]$, we have
\[
\|\eta_t-\eta^\star\|_{\mathrm{TV}}\, \underset{t \to \infty}{\longrightarrow}\, 0\, .
\]
\end{prop}

Before proving Proposition~\ref{prop:coupling-two-community}, let us state a useful lemma.

\begin{lemma}\label{lem:bound-(h+1)-gamma}
For $Y\sim\eta^\star$, provided $\gamma\leq \tfrac{1}{2}$, we have $\E[(|Y|+1)^\gamma]\leq 2^\gamma +\tfrac{\gamma}{\theta_1}$, with $\theta_1$ defined in~\eqref{eq:def-theta1}.
\end{lemma}

\begin{proof}[Proof]
Let $(Y_t)_{t\geq 0}$ be a solution of~\eqref{eq:new-dynamics-2comm} started at $Y_0\sim\eta^\star$. As in the proof of Proposition~\ref{prop:semiflow}, we start by observing that the height process $(|Y_t|)_{t\geq 0}$ is always stochastically dominated by the homogeneous Markov process $(H_t)_{t\geq 0}$ on $\N^*$ which starts from the height measure induced by $\eta^\star$, and which, when at $k=1$, jumps to $2$ at rate $2^\gamma$, and when at $k\geq 2$, jumps to $k+1$ at rate $(k+1)^\gamma$ and to $1$ at rate $\theta_1$. For this dominating process, we have
\[
\partial_t\E[(H_t+1)^\gamma]=\E\left[\left((H_t+2)^\gamma-(H_t+1)^\gamma\right)(H_t+1)^\gamma\right]-\theta_1\left(\E[(H_t+1)^\gamma]-2^\gamma\right)\, .
\]
Using $(H_t+2)^\gamma-(H_t+1)^\gamma\leq \gamma(H_t+1)^{\gamma-1}$, we have
\begin{align*}
\partial_t\E[(H_t+1)^\gamma] &\leq \gamma\E\left[(H_t+1)^{2\gamma-1}\right]-\theta_1\left(\E[(H_t+1)^\gamma]-2^\gamma\right)\\
&\leq -\theta_1\E[(H_t+1)^\gamma]+ 2^\gamma\theta_1+\gamma\, ,
\end{align*}
where for the last inequality, we used that $\gamma\leq \tfrac{1}{2}$. Hence $\E[(H_t+1)^\gamma]$ is upper-bounded by a function which converges to $2^\gamma+\tfrac{\gamma}{\theta_1}$. And since, for all $t\geq 0$, $\E[(|Y|+1)^\gamma]\leq \E[(H_t+1)^\gamma]$, we obtain the desired result.
\end{proof}

\begin{proof}[Proof of Proposition~\ref{prop:coupling-two-community}]
Let $(X_t, Y_t)_{t\geq 0}$ be a coupling of $(\eta_t)_{t\geq 0}$ and $\eta^\star$, generated as in Section~\ref{sec:thecoupling-2}, with $(X_0, Y_0)$ an arbitrary coupling of $\eta_0$ and $\eta^\star$ (\emph{e.g.} the independent coupling). We will show that for all $t>0$,
\[
\partial_t \P(X_t\neq Y_t)\leq -a\P(X_t\neq Y_t)\, ,
\]
with $a>0$ provided $\gamma$ is small enough. This will imply that $\P(X_t\neq Y_t) \underset{t\to \infty}{\longrightarrow} 0$, and hence that $\|\eta_t-\eta^\star\|_{\mathrm{TV}} \underset{t\to \infty}{\longrightarrow} 0$. 

In what follows, we use the notation
\[
\beta_\eta^\Lambda=\sum_{k\geq 1}\eta(k, \Lambda)\left((k+1)^\gamma-k^\gamma\right)=\alpha_\eta^\Lambda -p_\Lambda \, .
\]
One can already notice that, since the total-variation distance between two probability measures $\eta, \nu\in\cP(E)$ is the supremum, over all functions $f$ on $E$ taking values between $0$ and $1$, of the absolute difference $|\E_\eta[f]-\E_\nu[f]|$, and since the function $k\mapsto (k+1)^\gamma-k^\gamma$ is between $0$ and $2^\gamma -1$, we have 
\begin{equation}\label{eq:bound-beta-tv}
    |\beta_\eta^\Lambda -\beta_\nu^\Lambda|\leq (2^\gamma-1)\|\eta-\nu\|_{\mathrm{TV}}\, .
\end{equation}
Now, in order to control $\P(X_{t+\varepsilon}\neq Y_{t+\varepsilon}\given\cF_t)$, we observe the following:
\begin{itemize}
    \item when $X_t\neq Y_t$, then the two processes get together (either at $(1, A)$ or at $(1, B)$) at rate at least 
    \[
    \theta:= \sum_{\Lambda\in\{A, B\}}\inf_{\eta\in\cP(E)}s^\Lambda(\eta)\left(\alpha_\eta^\Lambda-3^\gamma +2^\gamma\right);
    \]
    \item when $X_t=Y_t=(1, \Lambda)$ with $\Lambda\in \{A, B\}$, then $X_{t+\varepsilon}\neq Y_{t+\varepsilon}$ when one of the two processes jumps upward to $(2, \Lambda)$ but not the other, which happens at rate $|s(\eta_t)-s(\eta^\star)|2^\gamma $;
    \item when $X_t=Y_t=(k, A)$ with $k\geq 2$, there are three main contributions to separation at time $t+\varepsilon$:
    \begin{itemize}
        \item either one jumps upward to $(k+1, A)$ but not the other, which happens at rate $|s(\eta_t)-s(\eta^\star)|(k+1)^\gamma $;
        \item or one jumps downward to $(1, A)$ but not the other, which happens at rate
        \begin{align*}
        &\left|s(\eta_t)(p_A+\beta_{\eta_t}^A-(k+1)^\gamma+k^\gamma)-s(\eta^\star)(p_A+\beta_{\eta^\star}^A-(k+1)^\gamma+k^\gamma)\right| \\
        &\hspace{1cm} \leq s(\eta_t)|\beta_{\eta_t}^A-\beta_{\eta^\star}^A|+|s(\eta_t)-s(\eta^\star)|(p_A+\beta_{\eta^\star}^A)\, ;
        \end{align*}
        \item or one jumps downward to $(1, B)$ but not the other, which happens at rate
        \begin{align*}
        &\left|(1-s(\eta_t))(p_B+\beta_{\eta_t}^B)-(1-s(\eta^\star))(p_B+\beta_{\eta^\star}^B)\right|\\
        &\hspace{1cm} \leq (1-s(\eta_t))|\beta_{\eta_t}^B-\beta_{\eta^\star}^B|+|s(\eta_t)-s(\eta^\star)|(p_B+\beta_{\eta^\star}^B)\, .
        \end{align*}
        \end{itemize}
        Summing those three possibilities, using the upper-bound~\eqref{eq:bound-beta-tv}, and the fact that $\beta_{\eta^\star}^A+\beta_{\eta^\star}^B\leq 2^\gamma-1$, we obtain a separation rate smaller than
        \begin{equation}\label{eq:separation-rate}
         |s(\eta_t)-s(\eta^\star)|\left((k+1)^\gamma +2^\gamma\right)+(2^\gamma-1)\|\eta_t-\eta^\star\|_{\mathrm{TV}}\, .
         \end{equation}
    \item when $X_t=Y_t=(k, B)$ with $k\geq 2$, one just has to adjust the preceding case but the same upper bound~\eqref{eq:separation-rate} holds for the separation rate.
\end{itemize}
All in all, we get that when the two processes are different at time $t$, they get together at rate at least $\theta$, whereas when they are together at $x\in E$, denoting by $|x|$ the height of $x$, they separate at rate at most 
\[
r_{|x|}(t):=(2^\gamma-1)\left(1+\frac{2^\gamma\left((|x|+1)^\gamma+2^\gamma\right)}{4p_Ap_B}\right)\|\eta_t-\eta^\star\|_{\mathrm{TV}}\, ,
\]
where we used the upper-bound~\eqref{eq:ineq-s(eta)-tv} on $|s(\eta_t)-s(\eta^\star)|$. For $\varepsilon > 0$, we then see that the probability that they are together at time $t+\varepsilon$ given that they are separated at time $t$ is larger than the probability that an exponential clock with parameter $\theta$ rings between $t$ and $t+\varepsilon$ (so that the two processes are together at level $1$), but another independent clock with parameter $2^\gamma$ (a crude upper-bound for the rate at which one of the two processes jumps from level $1$ to level $2$) does not ring in that same time interval:
\[
   \P(X_{t+\varepsilon} = Y_{t+\varepsilon}\given X_{t}\neq Y_{t})\geq (1-\e^{-\theta\varepsilon})\e^{-2^\gamma\varepsilon}=\theta\varepsilon+o(\varepsilon)\, ,
\]
as $\varepsilon \to 0$.
On the other hand, the probability that they are separated at time $t+\varepsilon$ given that they are together at $x$ at time $t$ is smaller than the probability that there exists a point in a zone with area $\int_t^{t+\varepsilon} r_{|x|}(s)\d s$:
\begin{align*}
   \P(X_{t+\varepsilon} \neq Y_{t+\varepsilon}\given X_{t}=Y_{t}=x) \leq 1-\e^{-\int_t^{t+\varepsilon} r_{|x|}(s)\d s}\leq \int_t^{t+\varepsilon} r_{|x|}(s)\d s\, .
\end{align*}
We thus have
\begin{align*}
    \P(X_{t+\varepsilon} \neq Y_{t+\varepsilon})&\leq \P(X_{t} \neq Y_{t})(1-\theta\varepsilon +o(\varepsilon))+\sum_{x\in E} \P(X_t=Y_t=x)\int_t^{t+\varepsilon}r_{|x|}(s)\d s \, .
\end{align*}
Using that $\P(X_t=Y_t=x)\leq \P(Y_t=x)=\eta^\star(x)$, we get
\begin{align*}
    \frac{\P(X_{t+\varepsilon} \neq Y_{t+\varepsilon})-\P(X_{t} \neq Y_{t})}{\varepsilon} \leq  -(\theta +o(1))\P(X_{t} \neq Y_{t})+\sum_{x\in E} \eta^\star(x)\frac{\int_t^{t+\varepsilon}r_{|x|}(s)\d s}{\varepsilon} \, .
\end{align*}
For $x\in E$ fixed, we have $\varepsilon^{-1}\int_t^{t+\varepsilon}r_{|x|}(s)\d s\to r_{|x|}(t)$ as $\varepsilon\to 0$. Moreover, 
\[
\varepsilon^{-1}\int_t^{t+\varepsilon}r_{|x|}(s)\d s \leq C({|x|}+1)^\gamma\, ,
\]
for some constant $C>0$, and Lemma~\ref{lem:bound-(h+1)-gamma} ensures that, if $\gamma\leq \tfrac{1}{2}$, then $\sum_{x\in E}\eta^\star(x)(|x|+1)^\gamma <+\infty$. Hence, by dominating convergence,
\[
\partial_t \P(X_t\neq Y_t)\leq -\theta\P(X_t\neq Y_t)+(2^\gamma-1)\left(1+\frac{2^\gamma\left(\E[(|Y|+1)^\gamma]+2^\gamma\right)}{4p_Ap_B}\right)\|\eta_t-\eta^\star\|_{\mathrm{TV}}\, ,
\]
where $Y\sim\eta^\star$. Finally, using Lemma~\ref{lem:bound-(h+1)-gamma} and the fact that $\|\eta_t-\eta^\star\|_{\mathrm{TV}}\leq \P(X_t\neq Y_t)$, we obtain
\[
\partial_t \P(X_t\neq Y_t)\leq -\left(\theta-(2^\gamma-1)\left(1+\frac{2^\gamma\left(2^{\gamma+1}+ \frac{\gamma}{\theta_1}\right)}{4p_Ap_B}\right)\right)\P(X_t\neq Y_t)\, ,
\]
and we see that we can find $\gamma$ small enough so that
\begin{equation}\label{condition-gamma-convergence}
\theta-(2^\gamma-1)\left(1+\frac{2^\gamma\left(2^{\gamma+1}+ \frac{\gamma}{\theta_1}\right)}{4p_Ap_B}\right)>0\, .
\end{equation}
For instance, for $p_A=p_B=\tfrac{1}{2}$ and $q_A=q_B=q$, taking $\gamma\leq \tfrac{1}{7}$ is sufficient. Indeed, in this case, we have $r_A=r_B$, and the dynamics is completely symmetric between $A$ and $B$. In particular, $s^A(\eta)=s^B(\eta)$ all along the trajectory. Using that $\alpha_\eta^A+\alpha_\eta^B\geq 1$, we then have
\[
\theta\geq \frac{1}{2}-3^\gamma+2^\gamma \quad \text{ and }\quad \theta_1\geq \frac{1}{2}\left(1-3^\gamma+2^\gamma\right)\, .
\]
One may then check that inequality~\eqref{condition-gamma-convergence} holds for $\gamma=1/7$.
\end{proof}

\subsection{Proof of Theorem~\ref{th:deg-distrib-2} and~\ref{th:deg-distrib-1}}\label{sec:proof-theorems}

In Proposition~\ref{prop:semiflow1}, we have shown that $(X(t))_{t\geq 0}$ is an asymptotic pseudo-trajectory of the semiflow induced by the well-posed dynamical system~\eqref{eq:dynamicsystem-1}. In particular, the set of limit points of $(X(t))_{t\geq 0}$ is included in the uniformly attracting set $\{\eta^\star\}$ of the semiflow generated by~\eqref{eq:dynamicsystem-1} (see Proposition~\ref{prop:coupling-one-community} for one community, Proposition~\ref{prop:coupling-two-community} for two communities). In particular, the only limit point of the compact sequence $(\eta_n)_{n\geq 0}$ is $\eta^\star$. For one community, we immediately have
\[
\frac{m_n}{n}=\frac{m_0}{n}+\eta_n R \underset{n\to +\infty}{\longrightarrow}\eta^\star R \quad \text{a.s.}
\]
For two communities, by Lemma~\ref{lem:m_n-eta_n}, we have
\[
\frac{m_n}{n}=\frac{m_0}{n}+\eta_n R+\frac{1}{n}\sum_{i=1}^n \left(B_i-p, p-B_i, 0, \dots\right)\underset{n\to +\infty}{\longrightarrow}\eta^\star R \quad \text{a.s.}\, ,
\]
by the law of large numbers.

\subsection{Expression of the limiting distribution $\mu^\star$}

We adopt the notation of Lemmas~\ref{lem:eigenmeasure} and~\ref{lem:fixedpoint} and recall the expression for the limit $\eta^\star$ obtained in Proposition~\ref{prop:coupling-two-community}. We have just shown in Section~\ref{sec:proof-theorems} that the limit $\mu^\star$ in Theorem~\ref{th:deg-distrib-2} is given by $\mu^\star = \eta^\star R$, which yields 
\[
\mu^\star(1, \Lambda)=p_\Lambda-\eta^*(1,\Lambda)=p_\Lambda\frac{\lambda^\star}{\lambda^\star+ s_\Lambda^\star}\, ,
\]
and for all $k \geq 2$,
\begin{align*}
    \mu^\star(k, \Lambda) &= \eta^\star(k-1,\Lambda)-\eta^\star(k,\Lambda)\\
    &= p_\Lambda \frac{\lambda^\star}{\lambda^\star+ s_\Lambda^\star k^\gamma}\prod_{j=1}^{k-1} \frac{s_\Lambda^\star j^\gamma}{\lambda^\star +s_\Lambda^\star j^\gamma}\,=\,  p_\Lambda \frac{\alpha_{\eta^\star}^\Lambda/r_\Lambda}{\alpha_{\eta^\star}^\Lambda/r_\Lambda +k^\gamma} \prod_{j=1}^{k-1}\frac{j^\gamma}{\alpha_{\eta^\star}^\Lambda/r_\Lambda +j^\gamma}\, ,
\end{align*}
where the last equality follows from Remark~\ref{rem:link-lambda-alpha}. In the single-community case, the limiting measure can be expressed more concisely as follows:
\[
\mu^\star(k) = \frac{\lambda^\star}{\lambda^\star + k^\gamma} \, \prod_{j=1}^{k-1} \frac{j^\gamma}{\lambda^\star + j^\gamma} = \frac{\alpha^\star}{\alpha^\star + k^\gamma} \, \prod_{j=1}^{k-1} \frac{j^\gamma}{\alpha^\star + j^\gamma}  \, ,
\]
for all $k \in \N^*$, since $\lambda^\star = \alpha_{\eta^\star} =: \alpha^\star$ in the one-community case.

\begin{remark}
In both the single- and two-community settings, the limiting distribution has a power-law tail if and only if $\gamma = 1$. Otherwise, it is stretched-exponential. This behavior is therefore specific to the case $\gamma = 1$. In particular, $\mu^\star(k, \Lambda)\propto k^{-3}$.
\end{remark}

\section{Model 2: Proof strategy and main results}\label{sec:model2}

In this section, we give all the elements needed to adapt the proof to \model\ $2$ in the case $q_A=q_B=q$. In this case, we still have
\[
\eta_{n+1}=\left(1-\frac{1}{n+1}\right)+\frac{1}{n+1}\delta_{X_{n+1}}\, ,
\]
but now we have
\[
\E[\delta_{X_{n+1}}\given \cF_n]=p_A\frac{m_n P^A}{m_n P^A(E)}+p_B\frac{m_n P^B}{m_n P^B(E)}\, ,
\]
where, for $\Lambda\in\{A, B\}$,
\[
P^\Lambda=(1-q)D^\Lambda +q D^{\overline{\Lambda}}\, .
\]
Setting $\tilde{Q}^\Lambda= RP^\Lambda$, and recalling that $\frac{m_n}{n}\approx \eta_n R$, we are then incited to write
\begin{equation}\label{eq:dynamic-2}
    \eta_{n+1}=\eta_n+\frac{1}{n+1}\left(p_A\frac{\eta_n \tilde{Q}^A}{\eta_n \tilde{Q}^A(E)}+p_B\frac{\eta_n \tilde{Q}^B}{\eta_n \tilde{Q}^B(E)}-\eta_n +\tilde{\varepsilon}_{n+1}\right)\, ,
\end{equation}
where $\tilde{\varepsilon}_{n+1}=\delta_{X_{n+1}}-p_A\frac{\eta_n \tilde{Q}^A}{\eta_n \tilde{Q}^A(E)}+p_B\frac{\eta_n \tilde{Q}^B}{\eta_n \tilde{Q}^B(E)}$. As for \model\ $1$, it is convenient to slightly change the step sizes and rewrite~\eqref{eq:dynamic-2} as
\[
\eta_{n+1}=\eta_n +\tilde{\gamma}_{n+1}\left\{G(\eta_n)+\tilde{\varepsilon}'_{n+1}\right\}\, ,
\]
where
\begin{eqnarray*}
\tilde{\gamma}_{n+1} &=& \frac{1}{n+1}\left(\frac{p_A}{\eta_n \tilde{Q}^A(E)}+\frac{p_B}{\eta_n \tilde{Q}^B(E)}\right)\, ,\\
\tilde{\varepsilon}'_{n+1}&=& \frac{\tilde{\varepsilon}_{n+1}}{\frac{p_A}{\eta_n \tilde{Q}^A(E)}+\frac{p_B}{\eta_n \tilde{Q}^B(E)}}\, ,
\end{eqnarray*}
and where
\begin{equation}\label{eq:def-G}
    G(\eta)= \tilde{s}(\eta)\left(\eta\tilde{Q}^A-\tilde{\alpha}_\eta^A \eta\right)
    +(1-\tilde{s}(\eta))\left(\eta\tilde{Q}^B-\tilde{\alpha}_\eta^B \eta\right)\, ,
\end{equation}
with $\tilde{\alpha}_\eta^\Lambda= \eta\tilde{Q}^\Lambda(E)$ and $\tilde{s}(\eta)=\tfrac{p_A \tilde{\alpha}_\eta^B}{p_A \tilde{\alpha}_\eta^B+p_B \tilde{\alpha}_\eta^A}$. Note that, since $\tilde{Q}^\Lambda= (1-q)Q^\Lambda+q Q^{\overline{\Lambda}}$, the function $G$ defined in~\eqref{eq:def-G} has the same form as $F$ defined in~\eqref{eq:def-F}. More precisely, we have
\begin{align*}
G(\eta) &=\left((1-q)\tilde{s}(\eta)+q(1-\tilde{s}(\eta))\right)\left(\eta Q^A-\alpha_\eta^A \eta\right)\\
&\hspace{1cm}+\left((q\tilde{s}(\eta)+(1-q)(1-\tilde{s}(\eta))\right)\left(\eta Q^B-\alpha_\eta^B \eta\right)\, ,
\end{align*}
so that one can go from $F$ to $G$ by merely replacing $s(\eta)$ by $(1-q)\tilde{s}(\eta)+q(1-\tilde{s}(\eta))$. 

In particular, considering the increasing sequence $(\tilde{\tau}_n)_{n \geq 0}$ defined by $\tilde{\tau}_0=0$ and for all $n\geq 1$,
\[
\tilde{\tau}_n = \tilde{\gamma}_1 + \tilde{\gamma}_2 + ... + \tilde{\gamma}_n\, ,
\]
and introducing the linearized version $(\tilde{X}(t))_{t\geq 0}$ of $(\eta_n)_{n \in \N}$ defined, for all $n \geq 1$ and all $t \in [\tilde{\tau}_n, \tilde{\tau}_{n+1}]$, by
\[
\tilde{X}(t) = \eta_n + \cfrac{t-\tilde{\tau}_n}{\tilde{\tau}_{n+1}-\tilde{\tau}_n}(\eta_{n+1}-\eta_n)\, ,
\]
one can easily adapt the arguments of Section~\ref{sec:stochastic-approximation} to show that the interpolated process $(\tilde{X}(t))_{t\geq 0}$ is an asymptotic pseudo-trajectory of the semiflow induced by $G$. 

Similarly, all the arguments of Section~\ref{sec:convergence-two-communities} apply to \model\ $2$, including~\eqref{eq:ineq-s(eta)-tv} which still holds. Indeed,
\begin{align*}
    &\left|(1-q)\tilde{s}(\eta)+q(1-\tilde{s}(\eta))-(1-q)\tilde{s}(\nu)-q(1-\tilde{s}(\nu))\right|\leq \left|\tilde{s}(\eta)-\tilde{s}(\nu)\right|\\
    &\hspace{2cm} \leq \frac{\left|p_A \tilde{\alpha}_\eta^B\left(p_A\tilde{\alpha}_\nu^B+p_B\tilde{\alpha}_\nu^A\right)-p_A\tilde{\alpha}_\nu^B\left(p_A\tilde{\alpha}_\eta^B+p_B\tilde{\alpha}_\eta^A\right)\right|}{\left(p_A\tilde{\alpha}_\eta^B+p_B\tilde{\alpha}_\eta^A\right)\left(p_A\tilde{\alpha}_\nu^B+p_B\tilde{\alpha}_\nu^A\right)}\\
    &\hspace{2cm} \leq \frac{p_Ap_B \left|\tilde{\alpha}_\nu^A\tilde{\alpha}_\eta^B-\tilde{\alpha}_\eta^A\tilde{\alpha}_\nu^B\right|}{\left(p_A\left((1-q)p_B+qp_A\right)+p_B\left((1-q)p_A+qp_B\right)\right)^2}\, .
\end{align*}
For the numerator, we have
\begin{align*}
    \left|\tilde{\alpha}_\nu^A\tilde{\alpha}_\eta^B-\tilde{\alpha}_\eta^A\tilde{\alpha}_\nu^B\right|&\leq \tilde{\alpha}_\nu^A\left|\tilde{\alpha}_\eta^B-\tilde{\alpha}_\nu^B\right|+\tilde{\alpha}_\nu^B\left|\tilde{\alpha}_\nu^A-\tilde{\alpha}_\eta^A\right|\\
    &\leq \left(\tilde{\alpha}_\nu^A+\tilde{\alpha}_\nu^B\right)(2^\gamma-1)\|\eta-\nu\|_{\mathrm{TV}}\\
    &\leq (2^\gamma-1)\|\eta-\nu\|_{1}\, ,
\end{align*}
where we used that $\left|\tilde{\alpha}_\eta^\Lambda-\tilde{\alpha}_\nu^\Lambda\right|\leq (2^\gamma-1)\|\eta-\nu\|_{\mathrm{TV}}$ and that 
\[
\tilde{\alpha}_\nu^A+\tilde{\alpha}_\nu^B=1+\sum_{x\in E}\nu(x)\left((|x|+1)^\gamma-|x|^\gamma\right)\leq 2^\gamma\leq 2\, .
\]
As for the denominator, it is equal to $\left(q(p_A^2+p_B^2)+2(1-q)p_Ap_B\right)^2$, which is larger $(2p_Ap_B)^2$, since $p_A^2+p_B^2\geq 2p_Ap_B$. All in all, we have
\[
\left|(1-q)\tilde{s}(\eta)+q(1-\tilde{s}(\eta))-(1-q)\tilde{s}(\nu)-q(1-\tilde{s}(\nu))\right|\leq \frac{2^\gamma-1}{4p_Ap_B}\|\eta-\nu\|_{1}\, .
\]

\section{Remaining proofs}\label{sec:proof}

\subsection{Proof of~\eqref{eq:ineq-s(eta)-tv} and~\eqref{eq:ineq-c(eta)-tv}}\label{sec:proof-ineq-s(eta)-c(eta)-tv}

We have
\begin{align*}
    |s(\eta)-s(\nu)|&=\frac{\left|r_A\alpha_\eta^B\left(r_A\alpha_\nu^B+r_B\alpha_\nu^A\right)-r_A\alpha_\nu^B\left(r_A\alpha_\eta^B+r_B\alpha_\eta^A\right)\right|}{\left(r_A\alpha_\eta^B+r_B\alpha_\eta^A\right)\left(r_A\alpha_\nu^B+r_B\alpha_\nu^A\right)}\\
    &=\frac{r_Ar_B\left|\alpha_\eta^B\alpha_\nu^A-\alpha_\nu^B\alpha_\eta^A \right|}{\left(r_A\alpha_\eta^B+r_B\alpha_\eta^A\right)\left(r_A\alpha_\nu^B+r_B\alpha_\nu^A\right)}\, \cdot 
\end{align*}
Since $\alpha_\eta^\Lambda$ and $\alpha_\nu^\Lambda$ are both larger than $p_\Lambda$, the denominator is larger than $(r_Ap_B+r_Bp_A)^2$, and we have
\[
\frac{r_Ar_B}{(r_Ap_B+r_Bp_A)^2}=\frac{1}{p_Ap_B}\frac{r_Ap_B\cdot r_Bp_A}{(r_Ap_B+r_Bp_A)^2}\leq \frac{1}{4p_Ap_B}\, \cdot 
\]
For the numerator, we see that
\begin{align*}
    \left|\alpha_\eta^B\alpha_\nu^A-\alpha_\nu^B\alpha_\eta^A \right|&\leq \alpha_\eta^B\left|\alpha_\nu^A-\alpha_\eta^A\right|+
    \alpha_\eta^A\left|\alpha_\eta^B-\alpha_\nu^B\right| \, .
\end{align*}
Since both $\left|\alpha_\nu^A-\alpha_\eta^A\right|$ and $\left|\alpha_\eta^B-\alpha_\nu^B\right|$ are smaller than $(2^\gamma-1)\|\eta-\nu\|_{\mathrm{TV}}$, and since $\alpha_\eta^A+\alpha_\eta^B\leq 1+2^\gamma-1=2^\gamma$, we have
\begin{align*}
    |s(\eta)-s(\nu)|&\leq \frac{2^\gamma(2^\gamma-1)\|\eta-\nu\|_{\mathrm{TV}}}{4p_Ap_B}= \frac{2^\gamma(2^\gamma-1)\|\eta-\nu\|_{1}}{8p_Ap_B}\, , 
\end{align*}
which establishes~\eqref{eq:ineq-s(eta)-tv}. Moving on to~\eqref{eq:ineq-c(eta)-tv}, we write $\lambda(\eta)=\alpha_\eta^As(\eta)+\alpha_\eta^B(1-s(\eta))$ and apply the triangle inequality to get
\begin{align*}
     |\lambda(\eta)-\lambda(\nu)|&\leq (\alpha_\eta^A+\alpha_\eta^B)|s(\eta)-s(\nu)|+s(\nu)|\alpha_\eta^A-\alpha_\nu^A|+(1-s(\nu))|\alpha_\eta^B-\alpha_\nu^B|\, .
\end{align*}
Combining the previous bound on $|s(\eta)-s(\nu)|$, the fact that $\alpha_\eta^A+\alpha_\eta^B\leq 2^\gamma$, and that both $\left|\alpha_\nu^A-\alpha_\eta^A\right|$ and $\left|\alpha_\eta^B-\alpha_\nu^B\right|$ are smaller than $(2^\gamma-1)\|\eta-\nu\|_{\mathrm{TV}}$, we obtain
\begin{align*}
     |\lambda(\eta)-\lambda(\nu)|&\leq (2^\gamma-1)\|\eta-\nu\|_{\mathrm{TV}}\left(1+\frac{2^{\gamma}}{4p_Ap_B}\right)\leq \frac{1}{2}\|\eta-\nu\|_{1}\left(1+\frac{1}{2p_Ap_B}\right)\, ,
\end{align*}
which establishes~\eqref{eq:ineq-c(eta)-tv}. \qed

\subsection{Proof of~\eqref{eq:convergence-noise}}\label{sec:proof-convergence-noise}

We start by decomposing
\begin{align*}
    \sup_{m\geq n}\left\| \sum_{i=n}^{m}\frac{1}{i+1}\varepsilon_{i+1}\right\|_1
    &\leq \underbrace{\sup_{m\geq n}\left\| \sum_{i=n}^m \frac{1}{i+1}(\delta_{X_{i+1}}-\E\left[\delta_{X_{i+1}}\given \cF_i\right])\right\|_1}_{(\texttt{I})}\\
    &\hspace{1cm}+ \underbrace{\max_{\Lambda\in\{A, B\}}  \sum_{i=n}^{+\infty}\frac{1}{i+1}\left\| \frac{\tilde{m}_i D^\Lambda}{\tilde{m}_i D^\Lambda(E)} -\frac{\eta_i Q^\Lambda}{\eta_i Q^\Lambda(E)}\right\|_1}_{(\texttt{II})}\, ,
\end{align*}
with $\tilde{m}_i=\tfrac{m_i}{i}$. For term $(\texttt{I})$, using the Burkholder--Davis--Gundy inequality, we see that there exists $C_2>0$ such that
\[
\E[(\texttt{I})^2]\leq C_2 \E\left[\sum_{i=n}^{+\infty}\frac{1}{(i+1)^2}\|\delta_{X_{i+1}}-\E\left[\delta_{X_{i+1}}\given \cF_i\right]\|_1^2\right]\, .
\]
Now we have
\begin{align*}
\E\left[\left\|\delta_{X_{i+1}}-\E\left[\delta_{X_{i+1}}\given \cF_i\right]\right\|_1^2\right] &=\E\left[\left(\sum_{x\in E}|\delta_{X_{i+1}}(x)-\E\left[\delta_{X_{i+1}}(x)\given \cF_i\right]|\right)^2\right]\leq 4\, ,
\end{align*}
where we used the triangle inequality and the fact that $\sum_{x\in E}\delta_{X_{i+1}}(x)=1$. Hence we get
\begin{equation}\label{eq:bound-I}
\E[(\texttt{I})^2]\leq 4C_2\sum_{i\geq n}\frac{1}{(i+1)^2}\leq \frac{4C_2}{n}\, \cdot 
\end{equation}
As for term $(\texttt{II})$, we have
\begin{align*}
    \sum_{i=n}^{+\infty}\frac{1}{i+1}\left\| \frac{\tilde{m}_i D^\Lambda}{\tilde{m}_i D^\Lambda(E)} -\frac{\eta_i Q^\Lambda}{\eta_i Q^\Lambda(E)}\right\|_1 &\leq \underbrace{\sum_{i=n}^{+\infty}\frac{1}{i+1}\left\| \frac{\tilde{m}_i D^\Lambda-\eta_i Q^\Lambda}{\eta_i Q^\Lambda(E)}\right\|_1}_{(\texttt{II.a})}\\ &\hspace{1cm}+\underbrace{\sum_{i=n}^{+\infty}\frac{1}{i+1}\|\tilde{m}_i D^\Lambda\|_1\left|\frac{1}{\tilde{m}_i D^\Lambda(E)}-\frac{1}{\eta_i Q^\Lambda(E)}\right|}_{(\texttt{II.b})}\, .
\end{align*}
For term $(\texttt{II.a})$, using the fact that $\eta_iQ^\Lambda(E)$ is lower bounded by $p_\Lambda$ and recalling that $\tilde{m}_iD^\Lambda-\eta_iQ^\Lambda=(\tilde{m}_i-\eta_i R)D^\Lambda$, we see that Lemma~\ref{lem:m_n-eta_n} together with a variance bound for the binomial distribution imply that
\begin{equation}\label{eq:bound-IIa}
\E\left[(\texttt{II.a})\right] \lesssim \sum_{i=n}^{+\infty} \frac{1}{(i+1)^{2}}\E\left[\big|\sum_{j=1}^i(B_j-p)\big|\right] \lesssim \sum_{i=n}^{+\infty} \frac{1}{(i+1)^{3/2}}\lesssim \frac{1}{\sqrt{n}}\, \cdot 
\end{equation}
As for term $(\texttt{II.b})$, noticing that $\|\tilde{m}_i D^\Lambda\|_1=\tilde{m}_i D^\Lambda(E)$, that $\eta_i Q^\Lambda(E)$ is larger than $p_\Lambda$, and applying Lemma~\ref{lem:m_n-eta_n} again, we have
\begin{equation}\label{eq:bound-IIb}
\E\left[(\texttt{II.b})\right] \lesssim \sum_{i=n}^{+\infty} \frac{1}{i+1}\E\left[|\eta_i Q^\Lambda(E)-\tilde{m}_i D^\Lambda(E)|\right]\lesssim \sum_{i=n}^{+\infty} \frac{1}{(i+1)^{3/2}}\lesssim \frac{1}{\sqrt{n}}\, \cdot 
\end{equation}
Combining \eqref{eq:bound-I} (together with the Jensen's inequality), \eqref{eq:bound-IIa} and \eqref{eq:bound-IIb}, we obtain
\[
 \E\left[ \sup_{m\geq n}\left\| \sum_{i=n}^{m}\frac{1}{i+1}\varepsilon_{i+1}\right\|_1\right]\lesssim \frac{1}{\sqrt{n}}\, \cdot 
\]
Markov's inequality then gives, for all $\varepsilon>0$,
\[
\sum_{k\geq 1}\P\left(\sup_{m\geq 2^k}\left\| \sum_{i=2^k}^{m}\frac{1}{i+1}\varepsilon_{i+1}\right\|_1>\varepsilon\right)\lesssim\sum_{k\geq 1}2^{-k/2}<+\infty\, ,
\]
so that, by the Borel--Cantelli lemma, we have $\sup_{m\geq 2^k}\left\| \sum_{i=2^k}^{m}\frac{1}{i+1}\varepsilon_{i+1}\right\|_1\to 0$ almost surely as $k\to +\infty$, and we can conclude by monotonicity. \qed

\subsection{Proof of Lemma~\ref{lem:fixedpoint}}\label{sec:proof-lem-fixedpoint}
We aim to show that the equation $s(\pi_u)=u$ admits a unique solution $u \in [0, 1]$. The following equivalences hold:
\begin{equation*}
\begin{aligned}
s(\pi_u) = u
&\iff r \pi_u Q^B(E)
= u(1-r) \pi_u Q^A(E) + ur \pi_u Q^B(E)\\
&\iff u \pi_u Q^A(E)
= r [
\underbrace{u \pi_u Q^A(E) + (1-u) \pi_u Q^B(E)}
_{=\pi_u Q_u(E) = \lambda(u)}]\\
&\iff \psi(u) := \cfrac{u \pi_u Q^A(E)}{u \pi_u Q^A(E) + (1-u) \pi_u Q^B(E)} = r \, .
\end{aligned}
\end{equation*}
Observe that the function $\psi$ is such that $\psi(0) = 0$ and $\psi(1) =1$. We will show that $\psi$ is strictly increasing and continuous by proving that its derivative is strictly positive. This will allow us to establish the existence and uniqueness of a solution $u$ to the equation $\psi(u) = r$ as $r \in [0, 1]$.\\
Let $\varphi$ be the function defined on $(0, \infty)$ by
\[
\varphi(x) = \displaystyle\sum\limits_{k=1}^{\infty} \displaystyle\prod\limits_{j=1}^{k} \cfrac{1}{1+x j^{-\gamma}}\, .
\]
Note that $\varphi$ is strictly decreasing and differentiable on $(0, \infty)$, with derivative given by
\[
\varphi'(x) = - \displaystyle\sum\limits_{k=1}^{\infty} \left(\displaystyle\prod\limits_{j=1}^{k} \cfrac{1}{1+x j^{-\gamma}}\right) \displaystyle\sum\limits_{i=1}^{k} \cfrac{1}{i^{\gamma}+x} \, .
\]
With this notation, a straightforward calculation shows that 
\[
\pi_u Q^A(E) = \cfrac{\lambda(u)}{u}p \, \varphi\left(\frac{\lambda(u)}{u}\right) \quad \text{and} \quad \pi_u Q^B(E) = \cfrac{\lambda(u)}{1-u}(1-p)\, \varphi\left(\frac{\lambda(u)}{1-u}\right) \, ,
\]
which yields the following simplified expression for $\psi$:
\[
\psi(u) = \cfrac{p \, \varphi\left(\frac{\lambda(u)}{u}\right)}{p \, \varphi\left(\frac{\lambda(u)}{u}\right) + (1-p) \, \varphi\left(\frac{\lambda(u)}{1-u}\right)} = p \, \varphi\left(\frac{\lambda(u)}{u}\right) \, ,
\]
since
\begin{align}\label{eq:def-lambda}
p \, \varphi\left(\frac{\lambda(u)}{u}\right) + (1-p) \, \varphi\left(\frac{\lambda(u)}{1-u}\right) = 1\, .
\end{align}
Now, let $h$ be the function defined by $h(u) = \frac{\lambda(u)}{u}$. We show that $h$ is strictly decreasing. For $u \in (0, 1)$, Equation~\eqref{eq:def-lambda} can be rewritten as
\[
p \, \varphi\left(h(u)\right) + (1-p) \, \varphi\left(h(u) \frac{u}{1-u}\right) = 1\, .
\]
Defining $f$ by
\[
\begin{array}{rcl}
f : (0, 1) \times (0, \infty) &\longrightarrow& \R \\
(u, y)              &\longmapsto& p \, \varphi(y) + (1-p) \, \varphi(y \times \frac{u}{1-u})-1 \, ,
\end{array}
\]
we have the following equivalence:
\[
f(u, y) = 0 \iff y = \frac{\lambda(u)}{u} \, .
\]
By the implicit function theorem, the derivative of $h$ is therefore given by
\[
h'(u) = - \cfrac{\cfrac{\partial f}{\partial u}(u, h(u))}{\cfrac{\partial f}{\partial y}(u, h(u))} 
= - \cfrac{\frac{(1-p)}{(1-u)^2}\varphi'\left(\frac{\lambda(u)}{1-u}\right)}{p \varphi'\left(\frac{\lambda(u)}{u}\right)+\frac{(1-p)u}{1-u}\varphi'\left(\frac{\lambda(u)}{1-u}\right)} \, .
\]
Since $h'$ is strictly negative on $(0, 1)$, $h$ is strictly decreasing. Moreover, as $\varphi$ is also strictly decreasing, it follows by composition that $\psi$ is strictly increasing. Standard arguments show that $\psi$ is continuous on $[0, 1]$. We have thus established the existence of a unique $s$ satisfying $\psi(u) = r$, and hence $s(\pi_u) = u$. This completes the proof. \qed

\bigskip

\subsection*{Acknowledgments} The authors warmly thank Pierre Monmarché for insightful discussions at various stages of this work and for his kindness. They also thank Nicolas Fournier for stimulating discussions on jump processes and for his interest in this work.

\bigskip

\bibliographystyle{abbrvnat}
\bibliography{references}

\end{document}